\documentclass[reqno,a4paper,11pt]{article}

\usepackage{rcs}
\RCS $Date: 2005/05/23 11:30:21 $
\RCS $Revision: 1.34 $

\usepackage{epic}
\usepackage{url}
\usepackage{hyperref}
\usepackage{amssymb,amsmath}
\usepackage{graphicx}
\usepackage{enumerate}
\usepackage[all]{xy}
\usepackage{color}
\usepackage{xcolor}
\usepackage{color}

\usepackage{amsthm}
\newtheorem{definition}{Definition}
\newtheorem{conjecture}[definition]{Conjecture}

\newtheorem{remark}[definition]{Remark}
\newtheorem{example}[definition]{Example}

\newtheorem{lemma}[definition]{Lemma}
\newtheorem{proposition}[definition]{Proposition}

\newtheorem{theorem}[definition]{Theorem}

\newcommand{\set}[1]{\left\{{#1}\right\}}

\newcommand{\vek}[1]{\boldsymbol{#1}}

\newcommand\setsuchas[2]{\left\{\,{#1}\,\vrule\,{#2}\,\right\}}

\newcommand{\Nat}{{\mathbb{N}}}
\newcommand{\Z}{{\mathbb{Z}}}
\newcommand{\C}{{\mathbb{C}}}
\newcommand{\R}{{\mathbb{R}}}

\newcommand{\CompOf}[2]{#1 \vDash #2}

\newcommand{\w}[1]{{\lvert #1 \rvert}}
\newcommand{\mw}{\mathbf{mw}}
\newcommand{\NC}{\mathfrak{N}}
\newcommand{\BBD}{\mathfrak{BBD}}
\newcommand{\STA}{\mathfrak{S}}
\newcommand{\UpSubs}{\Lambda}
\newcommand{\cover}{\gtrdot}
\newcommand{\Compositions}{\mathcal{C}}
\newcommand{\Young}{\ensuremath{\mathcal{Y}}}
\newcommand{\AlphabetLR}{\mathcal{A}_{LR}}
\newcommand{\AlphabetU}{\mathcal{A}_{U}}
\newcommand{\AlphabetV}{\mathcal{A}_{V}}
\newcommand{\Alphabet}{\mathcal{A}_{LRUV}}

\newcommand{\PBox}{\color{blue}\dashbox{1}(1,1)}
\newcommand{\PEBox}{\color{magenta}\dashbox{0.2}(1,1)}
\newcommand{\thisP}{      
  \multiput(0,0)(0,1.2){3}{\PBox}
  \multiput(1.2,0)(0,1.2){4}{\PBox}
  \multiput(2.4,0)(0,1.2){1}{\PBox}
  \multiput(3.6,0)(0,1.2){2}{\PBox}
}

\newcommand{\mybs}[1]{\boldsymbol{\color{magenta}#1}}
\newcommand{\Sort}{\mathfrak{K}}
\CompileMatrices

\begin{document}

\title{Saturated chains in composition posets}
\author{Jan Snellman \\
Department of Mathematics,
Stockholm University \\
SE-10691 Stockholm, Sweden \\
email: \texttt{Jan.Snellman@math.su.se}
}
\date{\RCSDate}

\maketitle

\begin{abstract}
  We study some poset structures on the set of all compositions. In
  the first case, the covering relation consists of inserting a part of
  size one to the left or to the right, or increasing the size of
  some part by one. The resulting poset \(\NC\) was studied by the
  author in \cite{Snellman:Ncterm} in relation to non-commutative term
  orders, and then in \cite{Snellman:StandardPaths}, where some
  results about generating functions for \emph{standard paths} in
  \(\NC\) was established. This was inspired by the work of Bergeron,
  Bousquet-M{\'e}lou and Dulucq \cite{StPa} on standard paths in the
  poset \(\BBD\), where there are additional cover relations which
  allows the insertion of a part of size one anywhere in the
  composition. Finally, following a suggestion by Richard Stanley we
  study a poset \(\STA\) which is an extension of \(\BBD\).
  This poset is related to
  quasi-symmetric functions.

  For  these posets, we study generating functions for saturated
  chains of fixed width \(k\). We also construct ``labeled''
  non-commutative generating 
  functions and their associated languages.
\end{abstract}




\begin{section}{Introduction}
  To an integer partition one can associate its \emph{diagram}, which is a
  finite subset of \(\Nat^2\). Ordering the set of partitions by
  inclusion of diagrams, one gets a locally finite, ranked, distributive
  lattice \(\Young\) which is known as \emph{Young's lattice}. The
  empty partition \(\emptyset\) is the unique minimal element, and
  saturated chains in \(\Young\) from the bottom element corresponds
  to an increasing sequence of diagrams, where at each step a single
  box is added. Such a sequence can be succinctly coded as a
  \emph{standard tableau} on the final diagram in the chain. This
  well-known construction is used not only in combinatorics, but also
  in the representation theory of the symmetric group.

  In \cite{Snellman:newglue} it was observed that Young's lattice
  also classifies the \emph{standard term orders}, i.e. admissible
  group orders \(\le\) on \(G=\Z^n\) where, for a fixed choice of basis
  \(\vek{e}_1,\dots,\vek{e}_n\) of \(G\) it holds that
  \(\vek{e}_1 < \dots < \vek{e}_n\). The correspondence is as
  follows: first, we can consider instead standard monoid orders on
  \(\Nat^n\), which is isomorphic to the monoid of power products in
  \(x_1,\dots,x_n\). Secondly, we send the power product \(x_i^r\) to
  the partition which has \(r\) parts of size \(i\). Third, we extend
  this to a monomial \(x_1^{a_1} \cdots x_n^{a_n}\) in the natural
  way. The image of this injective map will be all partitions with
  parts of size 
  \(\le n\). This is a sublattice of Young's lattice, and if we pull
  back this order to the monoid of power products in
  \(x_1,\dots,x_n\), we get a partial order which is the intersection
  of all standard term orders. Any standard term order is thus a
  \emph{multiplicative total extension} of this poset. We can be bold
  and allow infinitely many indeterminates \(x_1,x_2,x_3,\dots\): the
  poset so obtained is the isomorphic to  Young's lattice.

  If we do the same for non-commutative term orders, i.e. monoid
  orderings of the free non-commutative monoid on
  \(x_1,x_2,\dots\) such that \(x_1 < x_2 < \cdots\), then the
  resulting poset is no longer a 
  lattice. It is natural to map a non-commutative monomial to a
  \emph{composition} rather than a partition. If we order the set of
  compositions by pushing forward the order relation on
  non-commutative monomials via this bijection, then the
  resulting poset structure has the following covering relations:
  \begin{enumerate}
  \item \((1,a_1,\dots,a_n) \cover (a_1,\dots,a_n)\), i.e. we may
    insert a part of size one to the left,
  \item \((a_1,\dots,a_n,1) \cover (a_1,\dots,a_n)\), i.e. we may
    insert a part of size one to the right,
  \item   \((a_1,\dots,a_i,\dots,a_n) \cover
    (a_1,\dots,a_i+1,\dots,a_n)\), i.e. we may 
    increase the size of a  part by one.
  \end{enumerate}
  The ``sorting map'' from compositions to partitions is
  order-preserving, and we can regard the above poset as a
  non-commutative analogue of Young's lattice.

  It is, however, not the only possible such analogue! In \cite{StPa}, 
Bergeron,  Bousquet-M{\'e}lou and Dulucq consider an analogous poset
on the set of compositions.
This poset, henceforth denoted by \(\BBD\), is an extension of \(\NC\):
there are one additional type of
covering relations:
\begin{displaymath}
  (a_1,\dots,a_i,1,\dots,a_n,) \cover (a_1,\dots,a_n),
\end{displaymath}
i.e. one can insert a part of size 1 anywhere in the composition. They
encoded \emph{standard paths}, i.e. saturated chains from the empty
composition \(()\) to some composition \(P\) as tableau on the diagram
of \(P\). This is a direct counterpart to standard Young tableaux.

Using the theory of labeled binary trees, they were able to
explicitly solve the differential equation satisfied by the
exponential generating function for such standard paths, and
furthermore to give precise asymptotics for the number of such paths of
a given length.

They also considered the simpler problem of enumerating standard paths
of a fixed \emph{width} \(k\), i.e ending at a composition with \(k\)
parts. Here, the generating functions turned out to be rational, given
by a simple recurrence formula.

In \cite{Snellman:StandardPaths} the ideas of Bergeron 
\emph{et  al} were used to give generating functions for standard
paths of fixed width in \(\NC\). In the present paper, we consider
saturated chains starting from an arbitrary composition. We also
introduce a non-commutative generalization, which encodes all information
about the saturated chains, not only their endpoints. We'll see that
these non-commutative power series are still rational, hence
\emph{recognizable} and given by a finite state machine.

We also consider yet another poset, proposed by Richard Stanley.
This poset, which we denote by \(\STA^\infty\), extends \(\BBD\) in
such a way that a composition of \(n\) is covered by precisely \(n+1\)
compositions. It occurs naturally in the study of the fundamental
quasi-symmetric functions.
We consider an infinite family of posets \(\STA^d\), all
extending 
\(\BBD\), which have the desired poset \(\STA^\infty\) as their
inductive limit,
and introduce  a compact, unifying  formalism for describing these
posets, together with \(\NC\) and \(\BBD\).

Some related aspects of enumeration that we \emph{have not} addressed are 
\begin{itemize}
\item Enumeration of saturated chains without restriction of the
  width, as in \cite{StPa},
\item Non-saturated chains, i.e. with steps of length \(>1\), as in
  \cite{Gessel:Young,Stanley:Differential,Stanley:Differential2},
\item Oscillating tableaux, i.e. chains going up and down in th poset,
  as in \cite{Delest:Oscillating}.
\end{itemize}

\end{section}

\begin{section}{Posets of compositions}
  \begin{subsection}{Multi-rankings on compositions}

  By a \emph{composition} \(P\) we mean a sequence of positive
  integers 
  \begin{equation}
    \label{eq:P}
   P = (p_1,p_2,\dots,p_k),    
  \end{equation}
 which are the \emph{parts} of
  \(P\). We define the \emph{width} \(\ell(P)\) of \(P\) as the
  number of parts, and the \emph{height} as the size of the largest part.
  The \emph{weight} \(\w{P} = \sum_{i=1}^k p_k\) of \(P\)
  is the sum of its parts. If \(P\) has weight \(n\) then \(P\) is a
  composition of \(n\), and we write \(\CompOf{P}{n}\).
  
  Let \(\Compositions\) denote the set of all compositions
  (including the empty one). For a non-negative integer \(k\), let
  \(\Compositions_{(k)}\) denote the subset of compositions of width 
    \(k\).

  The \emph{diagram} of a composition
  \(P=(p_1,\dots,p_k)\) is the set of points \((i,j) \in \Z^2\)
  with \(1 \le j \le p_i\). Alternatively, we can replace the node
  \((i,j)\) by the square with corners
  \((i-1,j-1)\),\((i-1,j)\),\((i,j-1)\) and \((i,j)\). So the
  composition \((1,2,4)\) has diagram 
  \begin{center}
    \setlength{\unitlength}{0.4cm}
    \begin{picture}(4,5)
      \multiput(0,0)(0,1){1}{\PBox}
      \multiput(1,0)(0,1){2}{\PBox}
      \multiput(2,0)(0,1){4}{\PBox}
  \end{picture}.
  \end{center}

  Thus, for a composition \(P\) the \emph{height} and \emph{width} of
  \(P\) is the height and width of the smallest rectangle containing
  its diagram.

  \begin{definition}
  We let \(\Nat^\omega\) denote the poset of finitely supported maps
   \(\Nat^+ \to \Nat\), with component-wise comparison. The
   \emph{Young lattice}  \(\Young\) is the sublattice of (weakly)
   decreasing maps. For any positive integer \(n\), \(\Nat^n\) can be
   identified with the subposet of \(\Nat^\omega\) consisting of maps
   with support in \(\set{1,2,\dots,n}\).
  \end{definition}

    Let \(\vek{e}_i\) be the \(i\)'th unit vector, and
    put 
  \begin{equation}
    \label{eq:f}
    \vek{f}_j = \sum_{i=1}^j \vek{e}_i,
  \end{equation}
  We define the \emph{multi-weight} of \(P\) by 
  \begin{equation}
    \label{eq:mv}
    \mw(P) = \sum_{i=1}^k \vek{f}_{p_i} \in \Young
  \end{equation}
  We define \(\mw^*(P)=\mw(P)^*\), where \(*\) is  conjugation (which
  is an order-preserving involution on \(\Young\)).

  Note that  \(\mw^*(P)=(q_1,\dots,q_k)\) is the
  decreasing reordering of \(P=(p_1,\dots,p_k)\).

  \begin{definition}\label{def:multiranked}
    A locally finite poset \((A, \ge)\) is said to be
    \emph{\(\omega\)-multi-ranked} 
    if there exists a map 
    \begin{equation}
      \label{eq:multirankdef}
      \Phi: A \to \Nat^\omega
    \end{equation}
    such that 
    \begin{equation}
      \label{eq:mcovpfi}
      m \cover m' \quad \implies \quad \Phi(m) \cover \Phi(m') 
    \end{equation}

    The poset is \emph{\(n\)-multi-ranked} if it is
    \(\omega\)-multi-ranked by \(\Phi\) and \(\Phi(A) \subseteq \Nat^n\).
      
    Similarly, the poset is \(\Young\)-multi-ranked if
    \(\Phi(A) \subseteq \Young\).
  \end{definition}

  Clearly, if \(A\) is \(\omega\)-multiranked then it is ranked
  (i.e. 1-multiranked), 
  since \(\Phi\) followed by the collapsing
  \begin{equation}
    \begin{split}
      \Nat^\omega & \to \Nat \\
      \boldsymbol{\alpha} = (\alpha_1,\alpha_2,\dots) & \mapsto 
      \w{\boldsymbol{ \alpha}} = \sum_{i=1}^\infty \alpha_i
    \end{split}
  \end{equation}
  will be a ranking.

  We will presently introduce the partial orders \(\NC\) and \(\BBD\)
  on the set \(\Compositions\) of all compositions.
  Our main interest is the poset \(\NC\), which is related to
  non-commutative term orders \cite{Snellman:Ncterm}.
  This poset is \(\Young\)-ranked, as is the extension \(\BBD\),
  studied in \cite{StPa}.

  However, we
  will consider also an extension which we denote by \(\STA\). This
  poset, brought to our attention by Richard Stanley,  is not
  \(\Young\)-ranked. In order to have a concept broad 
  enough to also encompass  this poset, we define:    

\begin{definition}\label{def:almostmultiranked}
    A locally finite poset \((a, \ge)\) is said to be
    \emph{almost} \(\Young\)-multi-ranked
    if there exists 
    \begin{enumerate}
    \item an  extension \(T\) of the partial order \(\le\) on
      \(\Young\), and
    \item a map
    \begin{math}
      \label{eq:amultirankdef}
      \Phi: A \to (\Young,T)
    \end{math}
    such that
    \end{enumerate}
    \begin{enumerate}[i)]
    \item  \(T\) is  is rank-preserving, i.e. 
        \(\vek{\alpha} \cover  \vek{\beta}\) w.r.t. \(T\) implies that
        \(\w{\vek{\alpha}} = \w{\vek{\beta}} + 1\), 
    \item     
    \begin{math}
      \label{eq:acovpfi}
      m \cover m' \quad \implies \quad \Phi(m) \cover \Phi(m'),
    \end{math}
  \item \label{en:surj} 
    \(\Phi(A) = \Young\),
  \item \label{en:tight} 
    if \(\vek{\alpha} \cover  \vek{\beta}\) w.r.t. \(T\) then
    there are \(u,v \in A\) such that \(u \cover v\) and
    \(\Phi(u)=\vek{\alpha}\), 
    \(\Phi(v)=\vek{\beta}\). 
    \end{enumerate}
    
    In particular, \(A\) is ranked via \(u \mapsto \w{\Phi(u)}\).
  \end{definition}

  The definition of  \(\omega\)-ranked posets is very natural, the
  definition of \Young-ranked posets somewhat less so. The
  definition of almost \Young-ranked posets is very \emph{ad hoc}:
  it aims to capture enough of the salient features of the posets
  \(\STA^d\), to be defined below, so that they can be considered
  together with \(\NC\) and \(\BBD\). It is a matter of aesthetics if
  one includes the conditions \eqref{en:surj} and \eqref{en:tight}  or
  not. 
  \end{subsection}

\end{section}

\begin{section}{Some (almost)
    \protect\(\protect\Young\protect\)-multiranked posets of 
    compositions}
  We will define  posets \(\NC\), \(\BBD\), \(\STA^d\), with
  underlying set \(\Compositions\). 

  \begin{subsection}{Operations on compositions}

    We define the infinite alphabets 
    \begin{equation}
      \label{eq:3}
      \begin{split}
      \AlphabetLR &= \set{L,R}  \\
      \AlphabetU &= \setsuchas{U_j}{j \in \Nat^+} \\
      \AlphabetV^d &= \setsuchas{V_i^r}{i,r \in \Nat^+,\, r < d} \\
      \Alphabet & = \AlphabetLR \cup \AlphabetU \cup \AlphabetV^\infty
      \end{split}
    \end{equation}
    For an alphabet \(A\), we denote by \(A^*\) the free monoid on
    \(A\). 
    
  \begin{definition}
    We define the following partially defined operations on
    \(\Compositions\). Let \(P=(p_1,\dots,p_k)\) be a composition, then
    \begin{enumerate}
    \item \(L.P = (1,p_1,\dots,p_k)\), defined for all \(P\),
    \item \(R.P = (p_1,\dots,p_k,1)\), defined for all \(P\),
    \item \(U_j. P = (p_1,\dots, p_{j-1}, p_j +1, p_{j+1},
      \dots,p_k)\), defined when \(j \le k\),
    \item \(V_i^1 . P = (p_1,\dots,p_{i-1},1,p_i,\dots,p_k)\), defined
      when \(i \ge 2\), \(p_{i-1} \ge 2\),
    \item \(V_i^r . P = (p_1,\dots, p_{i-2}, p_{i-1} - r + 1, r,
      p_{i}, \dots, 
      p_k)\), 
      defined when \(i \ge 2\), \(p_{i-1}-r+1 \ge 2\), \(r \ge 2\).
    \end{enumerate}
    \end{definition}

    \begin{definition}
      We define  a partial left action on \(\Compositions\) by the
    free monoid \(\Alphabet^*\) in the following way.
    We define, for a word 
    \(w = w't\), where \(t \in \Alphabet\), \(P \in \Compositions\),
    \begin{equation}
      \label{eq:action}
      w.P = w'. (t.P)
    \end{equation}
    if the action of \(t\) on \(P\) is defined, and if
    recursively the 
    action of \(w'\) on \(t.P\) is defined.
  \end{definition}

  We give \(L\) and all the \(U_j\), the highest priority, followed 
  by the \(V_i^r\), with the convention that \(V_i^r\) has higher
  priority than \(V_j^s\) iff \(i < j\). The lowest priority is
  given to \(R\).

    \begin{definition}
    Let \(P\) be a composition and suppose that
    \begin{displaymath}
      \mathcal{B} \subseteq \Alphabet.
    \end{displaymath}
    \begin{enumerate}[i)]
    \item   The action of one of the operations above, call it
      \(T\), on a 
      composition \(P\) is \emph{admissible for \(P\)} (relative to
      \(\mathcal{B}\)) if it is defined, and if 
      \(T.P \neq S.P\) for all operations in \(\mathcal{B}\) with
      higher priority. 
    \item The action of a word \(W=VT \in \mathcal{B}^*\), \(T \in
      \mathcal{B}\), \(V \in \mathcal{B}^*\) is
      admissible for \(P\) (relative  \(\mathcal{B}\)) iff the action
      of \(T\) on \(P\) is admissible, and recursively the action of 
      \(V\) on \(T.P\) is admissible.
    \item   We let \(\langle \mathcal{B};\,\, P \rangle\) be the set
      of words in 
      \(\mathcal{B}^*\) that are admissible for \(P\). 
    \item   We let
      \(\le_{\mathcal{B}}\) be the smallest poset  
      \( \subset  \Compositions \times \Compositions\)
      which contains 
      \begin{displaymath}
        \setsuchas{(Q, w.Q)}{Q \in \Compositions, \, w \in \langle
          \mathcal{B}^*; Q \rangle} 
      \end{displaymath}
    \end{enumerate}
  \end{definition}

  So 
  \(P \le_{\mathcal{B}} Q\) if \(Q\) can be obtained from \(P\)
  using a sequence of admissible operations in \(\mathcal{B}\).

  \end{subsection}

  \begin{subsection}{A first example}
    \begin{lemma}
      The poset \((\Compositions, \le_{\AlphabetU})\)
      is isomorphic to the infinite direct sum 
      \begin{displaymath}
        \sum_{i \in \Nat} \Nat^i
      \end{displaymath}
    \end{lemma}
    \begin{proof}
      The map that sends the composition 
      \[\boldsymbol{\alpha}=(\alpha_1,\dots,\alpha_r)\] to 
      \[(\alpha_1-1,\dots,\alpha_r-1) \in \Nat^r\]
      is an order-preserving bijection.
    \end{proof}

    A part of the Hasse diagram of this non locallly finite poset is
    shown in Figure~\ref{fig:dirsum}. The other posets that we will
    introduce presently are all extensions of this posets, connecting
    the various components and also adding links within each component.
    \begin{figure}[t]
      \centering

  \setlength{\unitlength}{0.6cm}
      \begin{picture}(20,9)
        \put(1,1){\circle*{0.3}}
        \put(0.5,0.8){\(\emptyset\)}

        \multiput(4,1)(0,1){8}{\circle*{0.3}}
        \put(4,1){\line(0,1){5}}
        \dottedline{0.2}(4,6)(4,9)
        \put(3.5,0.8){1}
        \put(3.5,1.8){2}
        \put(3.5,2.8){3}
        \put(3.5,3.8){4}
        \put(3.5,4.8){5}
        \put(3.5,5.8){6}

        \multiput(11,1)(1,0){1}{\circle*{0.3}}
        \multiput(10,2)(2,0){2}{\circle*{0.3}}
        \multiput(9,3)(2,0){3}{\circle*{0.3}}
        \multiput(8,4)(2,0){4}{\circle*{0.3}}
        \multiput(7,5)(2,0){5}{\circle*{0.3}}
        \multiput(6,6)(2,0){6}{\circle*{0.3}}
        \put(11,1){\line(-1,1){6}}
        \put(12,2){\line(-1,1){5}}
        \put(13,3){\line(-1,1){4}}
        \put(14,4){\line(-1,1){3}}
        \put(15,5){\line(-1,1){2}}
        \put(16,6){\line(-1,1){1}}
        \put(11,1){\line(1,1){6}}
        \put(10,2){\line(1,1){5}}
        \put(9,3){\line(1,1){4}}
        \put(8,4){\line(1,1){3}}
        \put(7,5){\line(1,1){2}}
        \put(6,6){\line(1,1){1}}

        \put(10.3,0.6){11}
        \put(9.3,1.6){21}
        \put(12.2,1.6){12}
        \put(8.3,2.6){31}
        \put(10.2,2.7){22}
        \put(13.2,2.6){13}
      \end{picture}

      \caption{The poset \((\Compositions, \le_{\AlphabetU})\).}
      \label{fig:dirsum}
    \end{figure}

  \begin{subsubsection}{Graphical representations}
    We have already introduced the diagram of a composition. 
    Another graphical depiction is the so-called \emph{balls and bars}
    representation: here, the composition \(P=(p_1,\dots,p_r)\) is
    represented by \(r\) groups of balls, separated by vertical bars,
    the \(i\)'th group consisting of \(p_i\) balls.
    A third way of encoding the composition is to regard it as the
    ``index vector'' of a (non-commutative) monomial: the \(P\) above
    would be represented by 
    \begin{equation}
      \label{eq:ncmon}
    x_{p_1} \cdots x_{p_r} 
    \end{equation}

    The effect of the operations \(U_j\) on 
    \(P=(3,4,1,3)\)
    is as follows:

    \begin{center}
      \begin{tabular}{|c|c|l|l|l|}
        \hline
        Operation& Result          & Diagram & balls and bars &
        monomial  \\ \hline
      \(P\) & (3,4,1,2) &  
    \setlength{\unitlength}{0.25cm}
    \begin{picture}(5,5)
      \thisP
    \end{picture}
       & \(ooo|oooo|o|oo\) & \(x_3x_4x_1x_2\)
      \\        

      \(U_1.P\) & (4,4,1,2) & 
    \setlength{\unitlength}{0.25cm}
    \begin{picture}(5,5)
      \thisP
      \put(0,3.6){\PEBox}
    \end{picture}
      & \(\mybs{o}ooo|oooo|o|oo\)
      &\(x_4^2x_1x_2\) \\ 

      \(U_2.P\) & (3,5,1,2) & 
    \setlength{\unitlength}{0.25cm}
    \begin{picture}(5,6)
      \thisP
      \put(1.2,4.8){\PEBox}
    \end{picture}
      & \(ooo|\mybs{o}oooo|o|oo\)
      &\(x_3x_5x_1x_2\) \\ 

      \(U_3.P\) & (3,4,2,2) & 
    \setlength{\unitlength}{0.25cm}
    \begin{picture}(5,5)
      \thisP
      \put(2.4,1.2){\PEBox}
    \end{picture}
      & \(ooo|oooo|\mybs{o}o|oo\)
      &\(x_3x_4x_2^2\) \\ 

      \(U_4.P\) & (3,4,1,3) & 
    \setlength{\unitlength}{0.25cm}
    \begin{picture}(5,5)
      \thisP
      \put(3.6,2.4){\PEBox}
    \end{picture}
      & \(ooo|oooo|o|\mybs{o}o\)
      &\(x_3x_4x_1x_3\) \\ 

      \hline
      \end{tabular}
    \end{center}

    The operation \(U_j\) adds a
    box on top of the \(j\)'th column in the diagram, adds a ball to
    the \(j\)'th group of balls, and \textbf{replaces the \(i\)'th
      variable \(x_{p_j}\) in the monomial with the variable \(x_{p_j}+1\)}.
  \end{subsubsection}
    
  \end{subsection}

  \begin{subsection}{The posets \protect\(\protect\NC\protect\)}

  \begin{definition}
    We define the following poset on the underlying set
    \(\Compositions\) of compositions: 
    \begin{displaymath}
      \NC = (\Compositions, \le_{\mathcal{B}}),
    \end{displaymath}
    where 
    \begin{displaymath}
      \mathcal{B} = \left( \AlphabetLR \cup \AlphabetU \right).
    \end{displaymath}
  \end{definition}

  \begin{subsubsection}{Graphical representations of the operations
      \protect\(L\protect\) and \protect\(R\protect\)}

    The effect of the operations \(L,R\) on 
    \(P=(3,4,1,3)\)
    is as follows:

    \begin{center}
      \begin{tabular}{|c|c|l|l|l|}
        \hline
        Operation& Result          & Diagram & balls and bars &
        monomial  \\ \hline
      \(P\) & (3,4,1,2) &  
    \setlength{\unitlength}{0.25cm}
    \begin{picture}(6,5)(-0.5,0)
      \thisP
    \end{picture}
       & \(ooo|oooo|o|oo\) & \(x_3x_4x_1x_2\)
      \\        

      \(L.P\) & (1,3,4,1,2) & 
    \setlength{\unitlength}{0.25cm}
    \begin{picture}(6,5)(-0.5,0)
      \thisP
      \put(-1.2,0){\PEBox}
    \end{picture}
     & \(\mybs{o|}ooo|oooo|o|oo\) &
      \(x_1x_3x_4x_1x_2\) \\

      \(R.P\) & (3,4,1,2,1) & 
    \setlength{\unitlength}{0.25cm}
    \begin{picture}(6,5)(-0.5,0)
      \thisP
      \put(4.8,0){\PEBox}
    \end{picture}
     & \(ooo|oooo|o|oo\mybs{|o}\) &
      \(x_3x_4x_1x_2x_1\) \\

      \hline
      \end{tabular}
    \end{center}

    We see that \(L\) adds a box to the left of the diagram, inserts a
    \(o|\) to the left of the balls and bars, and multiplies the
    monomial to the left with \(x_1\). Similarly, 
    \(R\) adds a box to the right of the diagram, inserts a
    \(|o\) to the right of the boxes and bars, and multiplies the
    monomial to the right with \(x_1\). 
  \end{subsubsection}

  The poset \(\NC\) with covering relations given by the operations
  \(L,R,U_j\) was introduced in \cite{Snellman:Ncterm} as a poset on
  the free monoid \(X^*\), \(X=\set{x_1,x_2,x_3,\dots}\). It is the
  poset of all ``multiplicative consequences'' of the ordering
  \begin{displaymath}
    x_1 < x_2 < x_3 < x_4 < \cdots
  \end{displaymath}
  of the variables. For instance, 
  \begin{displaymath}
    x_2 < x_3 \quad \implies \quad x_1x_2x_5^2 < x_1x_3x_5^2 =
    U_2.(x_1x_2x_5^2). 
  \end{displaymath}
  Formally, it is the intersection of all standard term orders on
  \(X^*\), where a standard term order is a total order such that
  \begin{displaymath}
    p < q \quad \implies \quad upv < uqv, \qquad \forall p,q,u,v \in X^*.
  \end{displaymath}
  The beginning of the Hasse diagram of \(\NC\) is shown in
  Figure~\ref{fig:HasseNC4}. 
\end{subsection}

\begin{subsection}{The poset \protect\(\protect\BBD\protect\)}
  \begin{definition}
    We define
    \(\BBD (\Compositions, \le_{\mathcal{B}})\)
      where 
      \[\mathcal{B}= \left( \set{L} \cup \AlphabetU \cup \AlphabetV^2
      \right). \]  
    \end{definition}
      This is the poset studied in \cite{StPa}. Compared to \(\NC\),
      it has the additional covering relations given by \(V_i^1\) which
      inserts a part of size one after a part of size \(\ge
      2\). 
      Graphically, this looks    like

    \begin{center}
      \begin{tabular}{|c|c|l|l|l|}
        \hline
        Operation& Result          & Diagram & balls and bars &
        monomial  \\ \hline
      \(P\) & (3,4,1,2) &  
    \setlength{\unitlength}{0.25cm}
    \begin{picture}(5,5)
      \thisP
    \end{picture}
       & \(ooo|oooo|o|oo\) & \(x_3x_4x_1x_2\)
      \\        

      \(V_2^1.P\) & (3,1,4,1,2) & 
    \setlength{\unitlength}{0.25cm}
    \begin{picture}(5,5)
      \multiput(0,0)(0,1.2){3}{\PBox}

      \multiput(2.4,0)(0,1.2){4}{\PBox}
      \multiput(3.6,0)(0,1.2){1}{\PBox}
      \multiput(4.8,0)(0,1.2){2}{\PBox}
      \put(1.2,0){\PEBox}
    \end{picture}
     & \(ooo|\mybs{o|}oooo|o|oo\) &
      \(x_3x_1x_4x_1x_2\) \\

      \(V_3^1.P\) & (3,4,1,1,2) & 
    \setlength{\unitlength}{0.25cm}
    \begin{picture}(5,5)
      \multiput(0,0)(0,1.2){3}{\PBox}
      \multiput(1.2,0)(0,1.2){4}{\PBox}

      \multiput(3.6,0)(0,1.2){1}{\PBox}
      \multiput(4.8,0)(0,1.2){2}{\PBox}
      \put(2.4,0){\PEBox}
    \end{picture}
     & \(ooo|oooo|\mybs{o|}o|oo\) &
      \(x_3x_4x_1^2x_2\) \\

      \(V_5^1.P\) & (3,4,1,2,1) & 
    \setlength{\unitlength}{0.25cm}
    \begin{picture}(5,5)
      \multiput(0,0)(0,1.2){3}{\PBox}
      \multiput(1.2,0)(0,1.2){4}{\PBox}
      \multiput(2.4,0)(0,1.2){1}{\PBox}
      \multiput(3.6,0)(0,1.2){2}{\PBox}
      \put(4.8,0){\PEBox}
    \end{picture}
     & \(ooo|oooo|o|oo\mybs{|o}\) &
      \(x_3x_4x_1x_2x_1\) \\ \hline

    \end{tabular}
  \end{center}

  Note that \(V_4^1\) is not admissible for this \(P\), and that
  adding a part of size one to the right is represented by \(V_5^1\)
  rather than by \(R\); in general, if the composition has \(r\)
  parts, and ends with a run of \(k\)
  parts of size 1, adding a one to right is represented by
  \(V_{r-k+1}^1\).
\end{subsection}

\begin{subsection}{The posets \protect\(\protect\STA^d\protect\)}
  \begin{definition}
      For \(d\) a positive integer \(> 1\) or \(d=\infty\),
      \(\STA^d =(\Compositions, \le_{\mathcal{A}})\), where 
      \[\mathcal{A}= \set{L} 
      \cup \AlphabetU
      \cup \AlphabetV^d. \] 
  \end{definition}

  The operations \(V_j^2,V_j^3\) operate as follows on \(P=(3,4,1,2)\).
    \begin{center}
      \begin{tabular}{|c|c|l|l|l|}
        \hline
        Operation& Result          & Diagram & balls and bars &
        monomial  \\ \hline
      \(P\) & (3,4,1,2) &  
    \setlength{\unitlength}{0.25cm}
    \begin{picture}(5,5)
      \thisP
    \end{picture}
       & \(ooo|oooo|o|oo\) & \(x_3x_4x_1x_2\)
      \\        

      \(V_2^2.P\) & (2,2,4,1,2) & 
    \setlength{\unitlength}{0.25cm}
    \begin{picture}(5,5)
      \multiput(0,0)(0,1.2){2}{\PBox}

      \multiput(2.4,0)(0,1.2){4}{\PBox}
      \multiput(3.6,0)(0,1.2){1}{\PBox}
      \multiput(4.8,0)(0,1.2){2}{\PBox}

      \multiput(1.2,0)(0,1.2){2}{\PEBox}
    \end{picture}
     & \(oo\mybs{|oo}|oooo|o|oo\) &
      \(x_2^2x_4x_1x_2\) \\

      \(V_3^2.P\) & (3,3,2,1,2) & 
    \setlength{\unitlength}{0.25cm}
    \begin{picture}(5,5)
      \multiput(0,0)(0,1.2){3}{\PBox}
      \multiput(1.2,0)(0,1.2){3}{\PBox}
      \multiput(2.4,0)(0,1.2){2}{\PEBox}
      \multiput(3.6,0)(0,1.2){1}{\PBox}
      \multiput(4.8,0)(0,1.2){2}{\PBox}

    \end{picture}
     & \(ooo|ooo\mybs{|oo}|o|oo\) &
      \(x_3^2x_2x_1x_2\) \\

      \(V_3^3.P\) & (3,2,3,1,2) & 
    \setlength{\unitlength}{0.25cm}
    \begin{picture}(5,5)
      \multiput(0,0)(0,1.2){3}{\PBox}
      \multiput(1.2,0)(0,1.2){2}{\PBox}
      \multiput(2.4,0)(0,1.2){3}{\PEBox}
      \multiput(3.6,0)(0,1.2){1}{\PBox}
      \multiput(4.8,0)(0,1.2){2}{\PBox}
    \end{picture}
     & \(ooo|oo\mybs{|ooo}|o|oo\) &
      \(x_3x_2x_3x_1x_2\) \\

      \hline
      \end{tabular}
      \end{center}
      
      In contrast to the other operations, the \(V_i^r\)'s, with \(r
      \ge 2\), does not only involve adding an extra box to a column
      of the diagram, or inserting a new column; it also means
      \textbf{taking away} a box from the preceding column. This may seem
      unnatural and contrived, but there is another representation
      with respect to which these operations make perfect sense.
      
      A compositions \(\alpha=(\alpha_1,\dots,\alpha_r)\) of  \(n\) 
      can be encoded as a subset of \(\set{1,2,\dots,n-1}\) via the
      bijection
      \begin{equation}\label{eq:compbij}
        (\alpha_1,\dots,\alpha_r) \mapsto S_\alpha=\set{\alpha_1,\alpha_1 +
          \alpha_2, \dots, \alpha_1 + \cdots + \alpha_{r-1}}. 
      \end{equation}
      If \(\pi\) is a permutation on \(\set{1,2,\dots,n}\) which has
      \emph{descent set}  
      \begin{displaymath}
        D_\pi = \set{\alpha_1,\alpha_1 +
          \alpha_2, \dots, \alpha_1 + \cdots + \alpha_{r-1}},
      \end{displaymath}
      consider all permutations on \(\set{0,1,2,\dots,n}\) which can
      be obtained by inserting a zero anywhere in the one-line
      representation of \(\pi\). For each such permutation \(\tau\)
      (there are of course exactly \(n+1\) of them) calculate its
      descent set, and find the unique composition of \(n+1\) which
      maps to this descent set under \eqref{eq:compbij}. The
      compositions obtained are \emph{precisely} the compositions
      which cover \(\alpha\) in \(\STA\).

      \begin{example}
        Let \(P=(3,4,1,2)\) as before. This is represented as 
        \(S_P=\set{3,7,8} \subset \set{1,2,\dots,9}\). The permutation
        \(\pi=[1,2,4,3,5,6,9,8,7,10]\) has descent set \(S_P\). Inserting a
        zero at all possible places, we get 11 new permutations, 11
        new descent set, and finally 11 compositions covering \(P\),
        as shown in table~\ref{tab:stapo}
        \begin{table}[t]
          \centering
          \begin{tabular}{|llll|}
            \hline
            Permutation & Descent set & Composition & Operation\\
            \hline
            \(\lbrack0,1,2,4,3,5,6,9,8,7  , 10 \rbrack\) & 
            \(\set{4,8,9}\) & (4,4,1,2) &
            \(U_1.P\) \\
            
            \(\lbrack 1,0,2,4,3,5,6,9,8,7  , 10 \rbrack\) &
            \(\set{1,4,8,9}\) 
            & (1,3,4,1,2) & \(L.P\) \\
            
            \(\lbrack 1,2,0,4,3,5,6,9,8,7  , 10 \rbrack\) &
            \(\set{2,4,8,9}\) & (2,2,4,1,2) & \(V_2^2.P\)
            \\
            
            \(\lbrack 1,2,4,0,3,5,6,9,8,7  , 10 \rbrack\) &
            \(\set{3,8,9}\) & (3,5,1,2) & \(U_2.P\)
            \\
            
            \(\lbrack 1,2,4,3,0,5,6,9,8,7  , 10 \rbrack\) &
            \(\set{3,4,8,9}\) & (3,1,4,1,2) & \(V_2^1.P\)
            \\
            
            \(\lbrack 1,2,4,3,5,0,6,9,8,7  , 10 \rbrack\) &
            \(\set{3,5,8,9}\) & (3,2,3,1,2) & \(V_3^3.P\)
            \\
            
            \(\lbrack 1,2,4,3,5,6,0,9,8,7  , 10 \rbrack\) &
            \(\set{3,6,8,9}\) & (3,3,2,1,2) & \(V_3^2.P\)
            \\
            
            \(\lbrack 1,2,4,3,5,6,9,0,8,7  , 10 \rbrack\) &
            \(\set{3,7,9}\) & (3,4,2,2) & \(U_3.P\)
            \\
   
            \(\lbrack 1,2,4,3,5,6,9,8,0,7  , 10 \rbrack\) &
            \(\set{3,7,8}\) & (3,4,1,3) & \(U_4.P\)
            \\
            
            \(\lbrack 1,2,4,3,5,6,9,8,7,0  , 10 \rbrack\) &
            \(\set{3,7,8,9}\) & (3,4,1,1,2) & \(V_3^1.P\)
            \\

            \(\lbrack 1,2,4,3,5,6,9,8,7, 10, 0 \rbrack\) &
            \(\set{3,7,8,10}\) & (3,4,1,2,1) & \(V_4^1.P\)
            \\
            \hline
          \end{tabular}
          \caption{Cover of \((3,4,1,2)\) in \(\STA^\infty\)}
          \label{tab:stapo}
        \end{table}
      \end{example}
      
      It is implicit  in Stanley's book \cite{Stanley:En2} (see section
      7.19, and in particular exercise 7.93), that the \emph{fundamental
        quasi-symmetric functions} \(L_\alpha\) multiply according to 
      \begin{equation}
        \label{eq:qs}
        L_1 L_\alpha = \sum_{\beta \cover \alpha} L_\beta
      \end{equation}
      where \(\cover\) is the covering relation in \(\STA^\infty\).
      Here, the \(L_\alpha\)'s are defined by 
      \begin{equation}
        \label{eq:La}
        L_\alpha = \sum_{\substack{i_1 \le i_2 \le \cdots \le i_n \\
            i_j < i_{j+1} \text{ if } j \in S_\alpha}} x_{i_1} x_{i_2}
        \cdots x_{i_n}
      \end{equation}
      and the set 
      \begin{displaymath}
        \setsuchas{L_\alpha}{\CompOf{\alpha}{n}}
      \end{displaymath}
      is a basis for the homogeneous quasi-symmetric functions of
      degree \(n\).
\end{subsection}

\begin{subsection}{Multi-ranking}
  If we identify the  posets \(\NC\), \(\BBD\), and \(\STA^d\) with
  their graphs, which are subsets 
  of \(\Compositions \times  \Compositions\), then
  
  \begin{equation}
    \label{eq:composets}
    \NC \subsetneq  \BBD =  \STA^2 \subsetneq \STA^3
    \subsetneq \cdots \subsetneq  \cup_{d=1}^\infty \STA^d =\STA^\infty 
  \end{equation}

  The posets \(\NC\) and \(\BBD\) have the same Hasse diagram up to
  rank 4, shown in Figure~\ref{fig:HasseNC4}. For \(d>2\) there is an
  edge between \((3)\) and \((2,2)\) in \(\STA^d\).

  \setlength{\unitlength}{1cm}
  \begin{figure*}[t]
    \begin{center}
      \begin{picture}(12,6) 
        \put(5.9,0.7){\(()\)}
        \put(6,1){\line(0,1){1}}
        
        \put(5.3,2){\(1\)}
        \put(6,2){\line(-2,1){2}}
        \put(6,2){\line(2,1){2}}

        \put(3.2,2.6){\(11\)}
        \put(4,3){\line(-3,1){3}}
        \put(4,3){\line(1,1){1}}
        \put(4,3){\line(3,1){3}}

        \put(8.5,2.7){\(2\)}
        \put(8,3){\line(2,1){2}}
        \put(8,3){\line(-1,1){1}}
        \put(8,3){\line(-3,1){3}}

        \put(0.3, 4){\(111\)}
        \put(1,4){\line(-1,1){1}}
        \put(1,4){\line(1,1){1}}
        \put(1,4){\line(3,1){3}}
        \put(1,4){\line(5,1){5}}

        \put(4.2,4){\(12\)}
        \put(5,4){\line(-1,1){1}}
        \put(5,4){\line(1,1){1}}
        \put(5,4){\line(2,1){2}}
        \put(5,4){\line(3,1){3}}

        \put(6.3,4){\(21\)}
        \put(7,4){\line(-5,1){5}}
        \put(7,4){\line(-3,1){3}}
        \put(7,4){\line(0,1){1}}
        \put(7,4){\line(3,1){3}}

        \put(10.5,4){\(3\)}
        \put(10,4){\line(-2,1){2}}
        \put(10,4){\line(0,1){1}}
        \put(10,4){\line(2,1){2}}

        \put(0,5.3){\(1111\)}

        \put(2,5.3){\(211\)}

        \put(4,5.3){\(121\)}

        \put(6,5.3){\(122\)}
        
        \put(7,5.3){\(22\)}
        
        \put(8,5.3){\(13\)}

        \put(10,5.3){\(31\)}

        \put(12,5.3){\(4\)}
      \end{picture}

      \caption{The Hasse diagram of \(\NC\).}
      \label{fig:HasseNC4}
    \end{center}
  \end{figure*}

  \begin{lemma}[\cite{Snellman:Ncterm}]
    \(\NC\) and \(\BBD\) are \Young-multiranked.
  \end{lemma}
  \begin{proof}
    Let \(P=(p_1,\dots,p_k)\) be a composition, and put
    \[(q_1,\dots,q_k) = \mw^*(P),\]
    then \((q_1,\dots,q_k)\) is the
    decreasing reordering of \(P\). Adding a part of size 1 to \(P\)
    adds a part of size 1 at the end of \(\mw^*(P)\), and increasing a
    part by one increases one part of  \(\mw^*(P)\) by one (a part
    which is strictly greater than its right neighbor). These
    operations are covering relations in the Young lattice, and all
    covering relations can be achieved. Furthermore, \(\mw\) is
    surjective. 
  \end{proof}

  \begin{lemma}
      The posets \(\STA^d\) are almost \Young-multiranked.
  \end{lemma}
  \begin{proof}
    The operation \(\mw^*(P) \mapsto \mw^*(V_i^r.P)\) corresponds to
    \begin{multline*}
      (q_1,\dots , q_i, \dots,q_k) \mapsto
      (q_1,\dots,q_{i+1-r},r,\dots, q_k) \mapsto \\
      \mw^*(q_1,\dots,q_{i+1-r},r,\dots, q_k) \in \Young,
    \end{multline*}
    where the last step performs the necessary resorting so that the
    result is a partition. We can let \(\Young^d\) be the smallest
    poset containing the original relations of \(\Young\) together
    with these new ones. Then 
    \begin{displaymath}
      \mw^*: \Compositions \to \Young^d
    \end{displaymath}
    is an almost \(\Young\)-multiranking.
  \end{proof}

  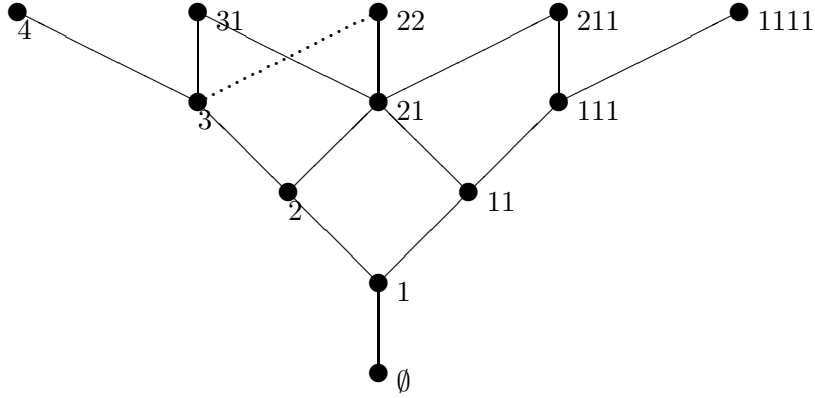
\begin{figure*}[t]
    \centering

  \setlength{\unitlength}{1.2cm}
    \begin{picture}(9,6)
      \put(4,1){\circle*{0.2}}
      \put(4,2){\circle*{0.2}}
      \multiput(3,3)(2,0){2}{\circle*{0.2}}
      \multiput(2,4)(2,0){3}{\circle*{0.2}}
      \multiput(0,5)(2,0){5}{\circle*{0.2}}
      
      \put(4.2,0.8){\(\emptyset\)}
      \put(4.2,1.8){1}
      \put(3.0,2.7){2}       
      \put(5.2,2.8){11}
      \put(2.0,3.7){3}       
      \put(4.2,3.8){21}       
      \put(6.2,3.8){111}
      \put(0,4.7){4}       
      \put(2.2,4.8){31}       
      \put(4.2,4.8){22}
      \put(6.2,4.8){211}
      \put(8.2,4.8){1111}

      \put(4,1){\line(0,1){1}}
      \put(4,2){\line(-1,1){2}}
      \put(2,4){\line(-2,1){2}}
      \put(4,2){\line(1,1){2}}
      \put(6,4){\line(2,1){2}}

      \put(3,3){\line(1,1){1}}
      \put(5,3){\line(-1,1){1}}

      \put(2,4){\line(0,1){1}}
      \dottedline[.]{0.1}(2,4)(4,5)
      \put(4,4){\line(-2,1){2}}
      \put(4,4){\line(0,1){1}}
      \put(4,4){\line(2,1){2}}

      \put(6,4){\line(0,1){1}}

    \end{picture}

    \caption{The lower part of the Hasse diagram of \(\Young^\infty\).
      The edges not present in \(\Young\) are  dotted.}
    \label{fig:YI}
  \end{figure*}

  \end{subsection}

  \begin{subsection}{Saturated chains, standard paths, and tableaux}
    Now suppose that \(\prec\) is a partial
    order on \(\Compositions\) such that  \(\mw\) is an
    almost \(\Young\)-multiranking on \((\Compositions, \prec)\).

    \begin{definition}
      If \(Q\) is a composition, we define a \emph{a saturated chain of length
        \(n\), starting from \(P\) and ending at \(Q\)} to be a
      sequence
      \begin{equation}
        \label{eq:gammachain}
      \gamma=(P_0=P,P_1,P_2,\dots,P_n=Q)
      \end{equation}
      of 
      compositions such that 
      \begin{equation}
        \label{eq:spath}
        P_0 \prec P_1 \prec P_2 \prec \cdots \prec P_n, \qquad \CompOf{P_i}{i},
      \end{equation}
      i.e. \(P_{i+1}\) should cover \(P_i\) for all \(i\).
      
      A \emph{standard path} is a saturated chain from the
      empty composition \(()\).
    \end{definition}

  We define the diagram, or the  \emph{shape}, of a
   a saturated chain from \(P\) to \(Q\) to be the diagram of \(Q\).
   
  Saturated chains in \(\NC\) or in \(\BBD\) can be coded
  as tableaux on the diagram of the terminal  composition.

  \begin{subsubsection}{\protect\(\protect\NC\protect\)}
  Consider first the poset \(\NC\).
  With respect to this order,
  \begin{equation} \label{eq:sp} 
  \rho=((1,2),(2,2),(1,2,2),(1,2,3),(1,2,4))
  \end{equation}
  is a saturated chain of length 
  4 from the minimal element 
  \((1,2)\) to the element \((1,2,4)\).

  It is clear what meant by saying that the boxes in the diagram of
  \(P_n\) should be labeled ``in the order that they 
  appear in the path'', if we
  (to avoid ambiguity)  use the convention that
  whenever \(P_i\) consists of \(i\) ones and \(P_{i+1}\) consists of
  \(i+1\) ones, the extra one is considered to have been added to the
  left. This in accordance with the above notion of priority of
  operations, since \(L\) has the highest priority.
  Thus, the only possible tableau  for standard paths ending at the
  composition \((1,1,1,1)\) is 
  \setlength{\unitlength}{0.35cm}
    \begin{picture}(4,1)
    \put(0,0){\line(1,0){4}}
    \put(0,1){\line(1,0){4}}
    \multiput(0,0)(1,0){5}{\line(0,1){1}}
    
    \put(0.2,0.2){\(4\)}
    \put(1.2,0.2){\(3\)}
    \put(2.2,0.2){\(2\)}
    \put(3.2,0.2){\(1\)}
  \end{picture}.

  As another example, the path \(\rho\) corresponds to the 
  tableau in Figure~\ref{fig:rhodiag}
  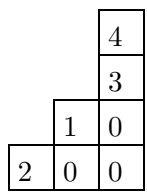
\begin{figure}[t]
    \centering
    \setlength{\unitlength}{0.6cm}
    \begin{picture}(3,5)
    \put(0,0){\line(1,0){3}}
    \put(0,1){\line(1,0){3}}
    \put(1,2){\line(1,0){2}}
    \put(2,3){\line(1,0){1}}
    \put(2,4){\line(1,0){1}}

    \put(0,0){\line(0,1){1}}
    \put(1,0){\line(0,1){2}}
    \put(2,0){\line(0,1){4}}
    \put(3,0){\line(0,1){4}}

    \put(1.2,0.2){0}
    \put(2.2,0.2){0}
    \put(2.2,1.2){0}
    \put(1.2,1.2){1}
    \put(0.2,0.2){2}
    \put(2.2,2.2){3}
    \put(2.2,3.2){4}
  \end{picture}    
    \caption{The diagram of the saturated chain \(\rho\) in \(\NC\).}
    \label{fig:rhodiag}
  \end{figure}
    
  \end{subsubsection}

  \begin{subsubsection}{\protect\(\protect\BBD\protect\)}
  Now consider the poset \(\BBD\), where there is the additional
  possibility of \(P < V_i^i.P\). If
  \(P=(p_1,p_2,\dots,p_\ell,1,\dots,1,p_s,\dots,p_n)\), with \(p_\ell
  >1\), then if \(\ell < i< j < s\) we 
have that \(V_i^1.P =
  V_j^1.P\). However, the  
  priority ordering and the rules for admissibility gives that only
  \(i=\ell+1\) is admissible. In other words, parts of size one can be
  inserted either to the left, or immediately after a part of size
  \(>1\). As an illustration, consider the following standard path
  (taken from \cite{StPa}):

  \begin{multline*}
    \gamma: \quad () \prec (1) \prec (1,1) \prec (1,2) \prec (1,2,1)
    \prec (1,3,1) \prec 
    (1,3,1,1) \prec \\ (1,3,1,2) \prec 
    (2,3,1,2) \prec (2,3,1,3)
    \prec (2,3,1,4) \prec (2,3,1,5) \prec (2,3,1,1,5)
  \end{multline*}

The diagram of \(\gamma\) is shown in Figure~\ref{fig:gammapath}.
\begin{figure}[t]
  \centering
    \setlength{\unitlength}{0.6cm}
    \begin{picture}(5,5)
    \multiput(0,0)(0,1){2}{\line(1,0){5}}
    \put(0,2){\line(1,0){2}}
    \put(1,3){\line(1,0){1}}
    \multiput(4,2)(0,1){4}{\line(1,0){1}}

    \put(0,0){\line(0,1){2}}
    \put(1,0){\line(0,1){3}}
    \put(2,0){\line(0,1){3}}
    \put(3,0){\line(0,1){1}}
    \put(4,0){\line(0,1){5}}
    \put(5,0){\line(0,1){5}}

    \put(0.2,0.2){2}
    \put(1.2,0.2){1}
    \put(2.2,0.2){12}
    \put(3.2,0.2){6}
    \put(4.2,0.2){4}
    \put(0.2,1.2){8}
    \put(1.2,1.2){3}
    \put(4.2,1.2){7}
    \put(1.2,2.2){5}
    \put(4.2,2.2){9}
    \put(4.2,3.2){10}
    \put(4.2,4.2){11}
  \end{picture}
  \caption{The diagram of the standard path \(\gamma\) in the poset \(\BBD\).}
  \label{fig:gammapath}
\end{figure}

  It is clear that for these two posets, given a tableau we can
  reconstruct the saturated chain.

  \emph{Nota bene}: the step   \(  (3,1,1,1) \cover (3,1,1)\) 
  is considered as adding a part at the right in \(\NC\), but as
  adding a part after \(3\) in \(\BBD\).

    
  \end{subsubsection}

  \begin{subsubsection}{The shadow of a tableau}
    \label{sec:shadow}
    Since the mapping \(\mw^*: \Compositions \to \Young\) is a
    multiranking for \(\NC\) and \(\BBD\), every saturated chain
    \(\gamma\), as in,  \eqref{eq:gammachain}
    in \(\NC\) or in \(\BBD\) ``lies over'' the saturated chain 
    \begin{equation}
      \label{eq:over}
      \mw^*(\gamma)=(\mw*(P_0),\mw^*(P_1),\mw^*(P_2),\dots,\mw^*(P_n))
    \end{equation}
    in \(\Young\). We call \(\mw^*(\gamma)\) the \emph{shadow} of
    \(\gamma\). If \(T\) is a tableau representing \(\gamma\), then we
    let \(\gamma(T)\) be the standard skew-tableau of shape
    \(\mw^*(P_n)/\mw^*(P_0)\) which encodes the way boxes are added to 
    \(\mw^*(P_0)\) to obtain \(\mw^*(P_n)\); we call this the shadow
    of \(T\). If \(S\) is a standard skew-tableau of shape
    \(\lambda/\mu\), then we define its multiplicity (w.r.t. \(\NC\)
    or \(\BBD\)) to be the number of saturated chains in \(\NC\) or in
    \(\BBD\) having \(\lambda/\mu\) as its shadow.
    
    \begin{example}
      The shadow of the tableaux in Figure~\ref{fig:rhodiag}
      is shown in Figure~\ref{fig:rhoshadow}. The shadow has
      multiplicity 8. The eight tableaux in \(\NC\) lying over the
      shadow is shown in Table~\ref{tab:allshad}.
      \begin{figure}[t]
        \centering
        \(S=\)
        \setlength{\unitlength}{0.6cm}
        \begin{picture}(3,5)
          \multiput(0,0)(0,1){4}{\PBox}
          \multiput(1,0)(0,1){2}{\PBox}
          \multiput(2,0)(0,1){1}{\PBox}
          \put(0.2,0.2){0}
          \put(0.2,1.2){0}
          \put(1.2,0.2){0}
          \put(1.2,1.2){1}
          \put(0.2,2.2){3}
          \put(0.2,3.2){4}
          \put(2.2,0.2){2}
        \end{picture}    
        \caption{The shadow of the the saturated chain \(\rho\) in
          \(\NC\).} 
        \label{fig:rhoshadow}
  \end{figure}

  \begin{table}[t]
    \centering

    \caption{The eight tableaux in \(\NC\) that lie over 
      the tableaux \(S\).}
      \begin{tabular}{llll}
      \setlength{\unitlength}{0.4cm}
      \begin{picture}(3,5)
        \multiput(0,0)(0,1){4}{\PBox}
        \multiput(1,0)(0,1){2}{\PBox}
        \multiput(2,0)(0,1){1}{\PBox}
        \put(0.2,0.2){0}
        \put(0.2,1.2){0}
        \put(1.2,0.2){0}
        \put(1.2,1.2){1}
        \put(0.2,2.2){3}
        \put(0.2,3.2){4}
        \put(2.2,0.2){2}
      \end{picture}    
      &
      \setlength{\unitlength}{0.4cm}
      \begin{picture}(3,5)
        \multiput(0,0)(0,1){2}{\PBox}
        \multiput(1,0)(0,1){4}{\PBox}
        \multiput(2,0)(0,1){1}{\PBox}
        \put(0.2,0.2){0}
        \put(0.2,1.2){0}
        \put(1.2,0.2){0}
        \put(1.2,1.2){1}
        \put(1.2,2.2){3}
        \put(1.2,3.2){4}
        \put(2.2,0.2){2}
      \end{picture}    
      &
      \setlength{\unitlength}{0.4cm}
      \begin{picture}(3,5)
        \multiput(0,0)(0,1){4}{\PBox}
        \multiput(1,0)(0,1){2}{\PBox}
        \multiput(2,0)(0,1){1}{\PBox}
        \put(0.2,0.2){0}
        \put(0.2,1.2){1}
        \put(1.2,0.2){0}
        \put(1.2,1.2){0}
        \put(0.2,2.2){3}
        \put(0.2,3.2){4}
        \put(2.2,0.2){2}
      \end{picture}    
      &
      \setlength{\unitlength}{0.4cm}
      \begin{picture}(3,5)
        \multiput(0,0)(0,1){2}{\PBox}
        \multiput(1,0)(0,1){4}{\PBox}
        \multiput(2,0)(0,1){1}{\PBox}
        \put(0.2,0.2){0}
        \put(0.2,1.2){1}
        \put(1.2,0.2){0}
        \put(1.2,1.2){0}
        \put(1.2,2.2){3}
        \put(1.2,3.2){4}
        \put(2.2,0.2){2}
      \end{picture}    
      \\
      \setlength{\unitlength}{0.4cm}
      \begin{picture}(3,5)
        \multiput(0,0)(0,1){1}{\PBox}
        \multiput(1,0)(0,1){2}{\PBox}
        \multiput(2,0)(0,1){4}{\PBox}
        \put(0.2,0.2){2}
        \put(1.2,0.2){0}
        \put(2.2,0.2){0}
        \put(1.2,1.2){1}
        \put(2.2,1.2){0}
        \put(2.2,2.2){3}
        \put(2.2,3.2){4}
      \end{picture}    
      &      
      \setlength{\unitlength}{0.4cm}
      \begin{picture}(3,5)
        \multiput(0,0)(0,1){1}{\PBox}
        \multiput(1,0)(0,1){4}{\PBox}
        \multiput(2,0)(0,1){2}{\PBox}
        \put(0.2,0.2){2}
        \put(1.2,0.2){0}
        \put(2.2,0.2){0}
        \put(1.2,1.2){1}
        \put(2.2,1.2){0}
        \put(1.2,2.2){3}
        \put(1.2,3.2){4}
      \end{picture}    
      &      
      \setlength{\unitlength}{0.4cm}
      \begin{picture}(3,5)
        \multiput(0,0)(0,1){1}{\PBox}
        \multiput(1,0)(0,1){2}{\PBox}
        \multiput(2,0)(0,1){4}{\PBox}
        \put(0.2,0.2){2}
        \put(1.2,0.2){0}
        \put(2.2,0.2){0}
        \put(1.2,1.2){0}
        \put(2.2,1.2){1}
        \put(2.2,2.2){3}
        \put(2.2,3.2){4}
      \end{picture}    
      &      
      \setlength{\unitlength}{0.4cm}
      \begin{picture}(3,5)
        \multiput(0,0)(0,1){1}{\PBox}
        \multiput(1,0)(0,1){4}{\PBox}
        \multiput(2,0)(0,1){2}{\PBox}
        \put(0.2,0.2){2}
        \put(1.2,0.2){0}
        \put(2.2,0.2){0}
        \put(1.2,1.2){0}
        \put(2.2,1.2){1}
        \put(1.2,2.2){3}
        \put(1.2,3.2){4}
      \end{picture}    
      \\      
      \end{tabular}
      \label{tab:allshad}
  \end{table}

    \end{example}
    
  \end{subsubsection}

  \end{subsection}
  
\end{section}

\begin{section}{Enumeration of saturated chains  of fixed width}

  \begin{subsection}{Definitions}

    Let \(\mathcal{Q}\) be one of the posets on compositions
    considered above.
    For a saturated chain \(\gamma\) of shape \((p_1,p_2,\dots,p_k)\) we
    set 
    \begin{equation}
      \label{eq:vwidth}
      v(\gamma)=x_1^{p_1} x_2^{p_2} \cdots x_k^{p_k}
    \end{equation}
    Note that this is a \emph{commutative} monomial, different
    from the representation used in \eqref{eq:ncmon}.

    We define the generating function
    \begin{equation}
      \label{eq:fk}
      f_k^{\alpha}(x_1,\dots,x_k) =
      f_k^{\alpha}[\mathcal{Q}](x_1,\dots,x_k) =
      \sum_{\substack{\gamma 
          \text{ saturated chain of width } k \\ 
          \text{ starting from } \alpha }}
      v(\gamma) 
    \end{equation}
    If \(\alpha\) is the empty composition, then we omit the
    superscript.

    \begin{definition}
      If \(f\) is a series in  \(x_1,x_2,x_3,\dots\) 
      and \(d,i,j\) are  positive integers, then
      \begin{equation}
        \label{eq:upsubs}
        \begin{split}
        \UpSubs_j(f)(x_1,x_2,x_3,\dots) &=
        f(x_1,x_2,\dots,x_{j-1},x_{j+1},x_{j+2},\dots) 
        \\
        \Delta_{i}^d(f) &= \frac{x_i^d}{d!} \frac{\partial_i^d f}{\partial
          x_i^d}(x_1,x_2,\dots,x_{i-1},0,x_{i+1},\dots)        
        \end{split}
      \end{equation}
    \end{definition}
     We put \(\UpSubs = \UpSubs_1\).

  \end{subsection}

  \begin{subsection}{Recurrence relations for the generating functions}

    \begin{lemma}\label{lemma:therecurs}
      Let \(P=(p_1,\dots,p_{k-1})\) be a composition. Then
      \begin{enumerate}[(i)]
      \item \(v(L.P) =x_1 \UpSubs(v(P))\),
      \item \(v(R.P) =x_k v(P)\),
      \item \(v(U_j.P) = x_j v(P)\), for \(j < k\),
      \item \(v(V_i^1.P) = x_i \UpSubs_i\bigl(v(P)- \Delta_{i-1}^1(v(P))
        \bigr) \), \(2 \le i \le k\),
      \item \(v(V_i^d.P) = \frac{x_i^d}{x_{i-1}^{d-1}}
        \UpSubs_i\bigl(v(P)- \sum_{j=1}^d \Delta_{i-1}^j(v(P)) 
        \bigr) \), \(2 \le i \le k\), \(d \ge 2\).
      \end{enumerate}
      if the respective operations are admissible (otherwise the RHS
      is zero).
    \end{lemma}
    \begin{proof}
      \(v(P) =x_1^{p_1} \cdots x_{k-1}^{p_{k-1}}\), so
      \begin{equation}
        \label{eq:aca}
        \begin{split}
        v(L.P)   &= v((1,p_1,\dots,p_{k-1})) \\
        &= x_1 x_2^{p_1} \cdots x_k^{p_{k-1}} \\
        &= x_1 \UpSubs(x_1^{p_1} \cdots x_{k-1}^{p_{k-1}}) \\
        &= x_1 \UpSubs(v(P)) \\ 
        v(R.P)   &= v((p_1,\dots,p_{k-1},1)) \\       
        &= x_1^{p_1} \cdots x_{k-1}^{p_{k-1}} x_k \\
        &= x_k v(P) \\ 
        v(U_j.P) &= v(p_1,\dots,p_{j-1}, p_j+1, p_{j+1},\dots,
        p_{k-1}) \\
        &= x_1^{p_1} \cdots x_{j-1}^{p_{j-1}} x_j^{p_j+1}
        x_{j+1}^{p_{j+1}} \cdots  x_{k-1}^{p_{k-1}} \\
        &= x_j x_1^{p_1} \cdots x_{j-1}^{p_{j-1}} x_j^{p_j}
        x_{j+1}^{p_{j+1}} \cdots  x_{k-1}^{p_{k-1}} \\
        &= x_j v(P) 
        \end{split}
      \end{equation}
      If \(p_{i-1} \ge 2\) then
      \begin{equation}
        \label{eq:tjo}
        \begin{split}
        v(V_i^1.P) &= v( p_1, \dots, p_{i-2}, p_{i-1} , 1, p_{i},
        \dots,  p_k) \\
        &= x_1^{p_1} \cdots x_{i-2}^{p_{i-2}} x_{i-1}^{p_{i-1}} x_i
        x_{i+1}^{p_{i}} \cdots    x_{k+1}^{p_k} \\ 
        &= x_i x_1^{p_1} \cdots x_{i-2}^{p_{i-2}} x_{i-1}^{p_{i-1}}
        \times x_{i+1}^{p_{i}} \cdots    x_{k+1}^{p_k} \\
        &= x_i \UpSubs_i(v(P)) \\
        &= x_i \UpSubs_i\left((v(P) - \Delta_{i-1}^1(v(P))\right) 
        \end{split}
      \end{equation}
      where the last equality follows from
      \begin{equation}
        \label{eq:tttjoo}
        \Delta_{i-1}^1(v(P)) = \Delta_{i-1}^1(m x_{i-1}^{p_{i-1}} m')
        = m 0 m' = 0
      \end{equation}
      On the other hand, if \(p_{i-1} <2\), i.e. if  \(p_{i-1} =1\),
      then
      \begin{equation}
        \label{eq:tttjoo2}
        \Delta_{i-1}^1(v(P)) = \Delta_{i-1}^1(m x_{i-1} m') = 
        x_{i-1} m m' = v(P)
      \end{equation}
      so
      \begin{equation}
        \label{eq:wewe}
        x_i \UpSubs_i\left((v(P) - \Delta_{i-1}^1(v(P))\right) = x_i
        \UpSubs_i\left(v(P) - v(P)\right) =0
      \end{equation}
      which is consistent with the fact that \(V_i^1\) is not
      admissible for \(P\).

      Similarly, to show that the action of \(V_i^d\) on \(P\)
      corresponds to 
      \[\frac{x_i^d}{x_{i-1}^{d-1}}
      \UpSubs_i\bigl(v(P)- \sum_{j=1}^d \Delta_{i-1}^j(v(P)) 
      \bigr)\] we want to show that
      \begin{equation}
        \label{eq:deltaip}
        \sum_{j=1}^d \Delta_{i-1}^j(v(P)) = 
         \begin{cases}
           0 & \text{ if } p_{i-1}-d+1 \ge 2 \\
           v(P) & \text{ otherwise }
         \end{cases}
      \end{equation}
      Suppose first that \(t=p_{i-1} \ge d+1\). Write \(v(P) =  m
      x_{i-1}^t\). 
      Then, for all \(1 \le j \le d\), 
      \begin{displaymath}
        \frac{\partial^j}{\partial x_{i-1}^j} m x_{i-1}^t = m j!
        x_{i-1}^{t-j} 
      \end{displaymath}
      is divisible by \(x_{i-1}\), hence 
      \begin{displaymath}
        \sum_{j=1}^d \Delta_{i-1}^j(v(P)) = 0.
      \end{displaymath}

      Suppose now that \(t=p_{i-1} \le d\). Then 
      \begin{equation}
        \label{eq:qwsss}
        \frac{\partial^j}{\partial x_{i-1}^j} m x_{i-1}^t = 
        \begin{cases}
          0 & j < t \\
          j! m & j= t \\
          0 & j > t 
        \end{cases}
      \end{equation}
      hence 
      \begin{equation}
        \label{eq:qwsss2}
        \Delta_{i-1}^j(m x_{i-1}^t ) =
        \begin{cases}
          0 & j < t \\
          m x_{i-1}^t & j= t \\
          0 & j > t 
        \end{cases}
      \end{equation}
      This proves the assertion.
    \end{proof}

    The above result gives recurrence relations for
    \(f_k^\alpha[\mathcal{Q}]\): 
    \begin{lemma}[\cite{Snellman:StandardPaths}]\label{lemma:ncrec}
      Let \(\mathcal{Q}=\NC\), let \(\alpha\) be a composition with
      \(r\) parts,
      and let \(f_k =  f_k^{\alpha}[\NC](x_1,\dots,x_k)\) be the generating
      function for saturated chains, starting from \(\alpha\), of width
      \(k\). Then \(f_k= 0\) for \(k < r\), and 
      \[f_r = v(\alpha) + (x_1 + \dots + x_r) f_r.\]
      Furthermore, for \(k > r\), \(f_k\) satisfies the following 
      recurrence relation
      \begin{equation}
        \label{eq:ncrec}
        f_k = x_1 \UpSubs(f_{k-1}) + x_k f_{k-1} + 
          (x_1 + \dots + x_k) f_k 
        \end{equation}
        if \(\alpha\) is not all-ones, and 
        \begin{equation}
          \label{eq:ncrec2}
          f_k =x_1 \UpSubs(f_{k-1}) + x_k f_{k-1} + (x_1 + \dots + x_k) f_k
          - x_1\cdots x_k 
        \end{equation}
         if \(\alpha\) is all-ones.
    \end{lemma}
    \begin{proof}
      This follows from Lemma~\ref{lemma:therecurs}, 
      since \(Q \in
      \Compositions_{(k)}\) can be obtained from \(P \in
      \Compositions_{(k-1)}\) either as \(Q=L.P\) or  
      \(Q=R.P\), and from \(W \in \Compositions_{(k)}\) as
      \(Q=U_i.W\). 
    \end{proof}

    \begin{lemma}[\cite{StPa}]
      \label{lemma:bbdrec}
    Let \(\mathcal{Q}=\BBD\),
    and let \(f_k =  f_k^{\alpha}[\BBD](x_1,\dots,x_k)\) be the generating
    function for saturated chains, from \(\alpha\) and of width
    \(k\). Then \(f_k= 0\) for \(k < r\), where \(r\) is the number of
    parts in \(\alpha\), and \[f_r = v(\alpha) + (x_1 + \dots + x_r) f_r.\]
    Furthermore, for \(k > r\), \(f_k\) satisfies the following
    recurrence relation
    \begin{equation}
      \label{eq:bbdrec}
      f_k =         x_1 \UpSubs(f_{k-1}) +  (x_1 + \dots + x_k) f_k  
        + \sum_{i=2}^k  x_i \UpSubs_i\left(f_{k-1}-
          \Delta_{i-1}^1(f_{k-1})\right) 
    \end{equation}
    \end{lemma}

    For the posets \(\STA^d\), the recurrence is as follows:
    \begin{lemma}
      \label{lemma:starec}
    Let \(\mathcal{Q}=\STA^d\), where \(d=\infty\) or \(d>2\) is a
    positive integer,
    and let \(f_k =  f_k^{\alpha}[\STA^d](x_1,\dots,x_k)\) be the generating
    function for saturated chains, from \(\alpha\) and of width
    \(k\). Then \(f_k= 0\) for \(k < r\), where \(r\) is the number of
    parts in \(\alpha\), and 
    \[f_r = v(\alpha) + (x_1 + \dots + x_r) f_r.\]
    Furthermore, for \(k > r\), \(f_k\) satisfies the following
    recurrence relation
    \begin{multline}
      \label{eq:starec}
      f_k = x_1 \UpSubs(f_{k-1}) +  (x_1 + \dots + x_k) f_k  
        + \\
        +
        \sum_{i=2}^k \sum_{1 \le v < d} 
        \frac{x_i^v}{x_{i-1}^{v-1}} \UpSubs_i \left(f_{k-1}-
          \sum_{j=1}^v \Delta_{i-1}^j(f_{k-1})  
        \right)
    \end{multline}
    \end{lemma}

    \begin{definition}
      Let \(\mathcal{Q}\) be one of the posets above, and
      let \(a_{n,k}^\alpha\) denote the number of saturated chains of
      width \(k\) and 
      length \(n-\w{\alpha}\), starting from \(\alpha\). 
      Define 
      \begin{equation}\label{eq:L}
        L_k^\alpha[\mathcal{Q}](t)=  L_k^\alpha(t) = \sum_{n\ge 0}
        a_{n,k}^\alpha t^n = f_k^\alpha(t,\dots,t) 
      \end{equation}
      \end{definition}

    Note that
    \(L_k^{()}=L_k^{(1)}\) for \(k >0\), so we may assume that \(\alpha\)
    has a positive number of parts.

\end{subsection}

\begin{subsection}{Enumeration of saturated chains of fixed width in
    the poset \protect\(\protect\NC\protect\)} 
The generating functions
\(f_k^\alpha[\NC] =f_k^\alpha\) are
displayed below
for some small \(k,\alpha\).  Note that \(f_k^{()} = f_k^{(1)}\)
for \(k > 0\). 

\begin{equation}
  \label{eq:fknc}
  \begin{split}
    f_1^{(1)} &= {\frac {x_{{1}}}{1-x_{{1}}}} \\
    f_2^{(1)} &= {\frac {x_{{1}}x_{{2}}\left (1-x_{{1}}x_{{2}}\right
    )}{\left (1-x_{{ 1}}\right )\left (1-x_{{2}}\right )\left
    (1-x_{{1}}-x_{{2}}\right )}} \\
f_2^{(1,1)} &= 
{\frac {x_{{1}}x_{{2}}}{1-x_{{1}}-x_{{2}}}} \\
f_3^{(1,1)} &= 
{\frac {x_{{1}}x_{{2}}x_{{3}}\left (1 -x_{{1}}x_{{2}} -x_{{1}}x_{{3}}
-{x_{{2}}}^{2} - x_{{2}}x_{{3}}\right )}{\left (1 - x_{{1}} -x_{{2}}\right )
\left (1 - x_{{2}} - x_{{3}}\right )
\left (1 - x_{{1}} -x_{{2}} -x_{{3}}\right )}} \\
f_3^{(1,2)} &= 
{\frac {x_{{1}}x_{{2}}x_{{3}}
\left ({-x_{{2}}}^{2} -2\,x_{{2}}x_{{3}}
+x_{{2}} -x_{{1}}x_{{3}} + x_{{3}}\right )}{\left (1- x_{{1}} -x_{{2}}\right )
\left (1 -x_{{2}} -x_{{3}}\right )\left (1 -x_{{1}} -x_{{2}} -x_{{3}}\right )}}
  \end{split}
\end{equation}

Using  the recurrence relation \eqref{eq:ncrec} we can prove by induction
\begin{theorem}[\cite{Snellman:StandardPaths}]\label{thm:fkform}
  For each \(k\),
  \begin{equation}
  \label{eq:fkstructureform}
  f_k^{()}(x_1,\dots,x_k) = \frac{x_1\cdots x_k }{
    \prod_{i=1}^k \prod_{j=i}^k (1-x_i -x_{i+1} - \ldots -x_{j})
    } \tilde{f}_k(x_1,\dots,x_k)
\end{equation}
where \(\tilde{f}_k\) is a polynomial.  
\end{theorem}

The corresponding result for \(f_k^{\alpha}\) is as follows:
\begin{theorem}\label{thm:fkalpha}
  Let \(\alpha\) be a composition with \(r>0\) parts. 
  For each \(k \ge r\),
  \begin{equation}
  \label{eq:fkstructurealpha}
  f_k^{\alpha}(x_1,\dots,x_k) = \frac{x_1\cdots x_k }{
    \prod_{i=1}^k \prod_{j=i+r-1}^k (1-x_i -x_{i+1} - \ldots -x_{j})
    } \tilde{f}_k^{\alpha}(x_1,\dots,x_k)
\end{equation}
where \(\tilde{f}_k^{\alpha}\) is a polynomial.  
\end{theorem}

\begin{proof}
  When \(k=r\) we have that \(f_k^\alpha = v(\alpha)(1-x_1-x_2- \cdots
  - x_k)^{-1}\), 
which has the desired form. For \(k > r\), assume 
  that \(f_k^\alpha\) has the above form. 
  
  If  \(\alpha\) is  all-ones, then
  by the recurrence relation
  \eqref{eq:ncrec2} it follows that
  \begin{multline}
      f_k^{\alpha} (1-x_1 - \cdots -x_k) =  
          x_1f_{k-1}^{\alpha}(x_2,\dots,x_k) +
          x_kf_{k-1}^{\alpha}(x_1,\dots,x_{k-1}) - x_1\cdots x_k  \\
          = x_1 x_2 \cdots x_k \tilde{f}_{k-1}^{\alpha}(x_2,\dots,x_k) 
          \prod_{i=2}^k \prod_{j=i+r-1}^k (1-x_i-\cdots -x_j)^{-1}
          \\
           +x_k x_1 \cdots x_{k-1} \tilde{f}_{k-1}^{\alpha}(x_1,\dots,x_{k-1}) 
          \prod_{i=1}^{k-1} \prod_{j=i+r-1}^{k-1} (1-x_i-\cdots -x_j)^{-1}
          - x_1\cdots x_k \\
          = x_1\cdots x_k \Biggl[ 
          \tilde{f}_{k-1}^{\alpha}(x_2,\dots,x_k) 
          \prod_{i=2}^k \prod_{j=i+r-1}^k (1-x_i-\cdots -x_j)^{-1}
          + \\
          \tilde{f}_{k-1}^{\alpha}(x_1,\dots,x_{k-1}) 
          \prod_{i=1}^{k-1} \prod_{j=i+r-1}^{k-1}
          (1-x_i-\cdots -x_j)^{-1}
          -1
          \Biggr] 
  \end{multline}
  hence 
  \begin{multline}\label{eq:multirec}
\tilde{f}_k^\alpha =    
\frac{f_k (1-x_1 - \cdots -x_k) 
\prod_{i=1}^k \prod_{j=i+r-1}^k
(1-x_i - \cdots -x_j)}
{x_1 \cdots x_k} 
  \end{multline}
is a polynomial.

If \(\alpha = (1,\dots,1)\) is
all-ones and has \(r\) parts, then
  by the recurrence relation
  \eqref{eq:ncrec} it follows that
  \begin{multline}
      f_k^{\alpha} (1-x_1 - \cdots -x_k) =  
          x_1f_{k-1}^{\alpha}(x_2,\dots,x_k) +
          x_kf_{k-1}^{\alpha}(x_1,\dots,x_{k-1})   \\
          = x_1 x_2 \cdots x_k \tilde{f}_{k-1}^{\alpha}(x_2,\dots,x_k) 
          \prod_{i=2}^k \prod_{j=i+r-1}^k (1-x_i-\cdots -x_j)^{-1}
          \\
           +x_k x_1 \cdots x_{k-1} \tilde{f}_{k-1}^{\alpha}(x_1,\dots,x_{k-1}) 
          \prod_{i=1}^{k-1} \prod_{j=i+r-1}^{k-1} (1-x_i-\cdots -x_j)^{-1}
          \\
          = x_1\cdots x_k \Biggl[ 
          \tilde{f}_{k-1}^{\alpha}(x_2,\dots,x_k) 
          \prod_{i=2}^k \prod_{j=i+r-1}^k (1-x_i-\cdots -x_j)^{-1}
          + \\
          \tilde{f}_{k-1}^{\alpha}(x_1,\dots,x_{k-1}) 
          \prod_{i=1}^{k-1} \prod_{j=i+r-1}^{k-1}
          (1-x_i-\cdots -x_j)^{-1}
          \Biggr] 
  \end{multline}
  hence 
  \begin{multline}\label{eq:multirec2}
\tilde{f}_k^\alpha =    
\frac{f_k (1-x_1 - \cdots -x_k) 
\prod_{i=1}^k \prod_{j=i+r-1}^k
(1-x_i - \cdots -x_j)}
{x_1 \cdots x_k} 
  \end{multline}
is a polynomial. 

\end{proof}

The generating functions
\[L_k^\alpha[\NC](t) = L_k^\alpha(t) = F_k^\alpha(t,\dots,t)\]
are clerly rational functions. 
We have that
\begin{equation}
  \label{eq:lksmall}
  \begin{split}
    L_1^{(1)} &= {\frac {t}{1-t}} \\
    L_2^{(1)} &= {\frac {\left (t+1\right ){t}^{2}}{\left
          (1-2\,t\right ) \left (1-t \right )}} \\ 
    L_3^{(1)} &= {\frac {{t}^{3}\left (-2\,{t}^{2}+5\,t+1\right
        )}{\left ( 1- 3\,t\right ) \left (1-2\,t\right )
        \left   (1-t\right )}} \\  
L_4^{(1)} &= {\frac {\left (6\,{t}^{3}-15\,{t}^{2}+16\,t+1\right
    ){t}^{4}}{\left (1-4 \,t\right )
\left (1-3\,t\right )\left (1-2\,t\right )
\left (1-t\right )}} \\
L_5^{(1,1)} &= {\frac {\left (-24\,{t}^{3}+38\,{t}^{2}-27\,t-1\right
    ){t}^{5}}{\left ( 1-5\,t\right )
\left (1-4\,t\right )
\left (1-3\,t\right )
\left (1-2\,t\right )}} \\
L_5^{(2,3)} &= 
8\,{\frac {{t}^{8}}{\left (1-5\,t\right )\left (1-4\,t\right )
\left (1-3\,t\right )\left (1-2\,t\right )}}
  \end{split}
\end{equation}

\begin{lemma}
  Let \(r\) denote the  number of parts
    of \(\alpha\), \(\CompOf{\alpha}{N}\).
  Then the 
  following recurrence relation holds:
  \begin{equation}
    \label{eq:fkrel2}
    \begin{split}
        L_r^\alpha &= t^N(1-rt)^{-1} \\
        L_k^\alpha &= \frac{2tL_{k-1}^\alpha -t^k}{1-kt}, \qquad k >
        r, \, \alpha \text{ all-ones} \\
        L_k^\alpha &= \frac{2tL_{k-1}^\alpha }{1-kt}, \qquad k >
        r, \, \alpha \text{ not all-ones} 
    \end{split}
  \end{equation}
\end{lemma}
\begin{proof}
  Specialize \eqref{eq:ncrec} and \eqref{eq:ncrec2}.
\end{proof}

We get  by induction:

\begin{lemma}
  Suppose that \(\alpha\) is not all-ones.
  Then 
  \begin{equation}
    \label{eq:10}
    L_k^\alpha = 2^{k-r}t^{N+k-r}\prod_{i=r}^k(1-it)^{-1}
  \end{equation}
\end{lemma}

Since 
\begin{equation}
  \label{eq:2}
  \prod_{j=r}^k (1-jt)^{-1} = \frac{k^{k-r}}{(k-r)!} (1-kt)^{-1} + l.o.t,
\end{equation}
we get that, when \(\alpha\) is a composition of \(N\) with \(r < N\)
parts,
\begin{equation}
  \label{eq:anka}
  a_{n,k}^{\alpha} \sim \frac{2^{k-r}}{k^N (k-r)!}
 k^n 
  \qquad \text{ as } n \to \infty 
\end{equation}

Now suppose that \(\alpha\) is all-ones, i.e. \(N=r\).
\begin{proposition}
  Suppose that \(\alpha\) is all-ones and has \(r>0\) parts. Then
\begin{equation}
  \label{eq:Lksa}
  L_k^\alpha(t) = \frac{t^k D_k^\alpha(t)}{\prod_{i=r}^k(1-it)} 
\end{equation}
where \(D_k^\alpha(t)\) is a polynomial satisfying the
recurrence
\begin{equation}
  \label{eq:11}
  D_k^\alpha(t) = 2 D_{k-1}^\alpha(t) - \prod_{i=r}^{k-1} (1-it).
\end{equation}
with initial conditions \(D_r^\alpha(t)=1\).
\end{proposition}
\begin{proof}
  This is true for \(k=r\). The assertion follows by induction, the
  induction step being
  \begin{equation}
    \begin{split}
      L_k &= \frac{t^k D_k}{\prod_{i=r}^k(1-it)} \\
          &= \frac{2t L_{k-1} - t^k}{1-kt} \\
          &= \frac{2t \frac{t^{k-1} D_{k-1}}{\prod_{i=r}^{k-1}(1-it)}
            -t^k}{1-kt} \\
          & = \frac{2t^k D_{k-1} - t^k \prod_{i=r}^{k-1}
            (1-it)}{\prod_{i=r}^{k} (1-it)}  
    \end{split}
  \end{equation}
  from which \eqref{eq:11} follows.
\end{proof}

The following proposition is a generalization of
a result in \cite{Snellman:StandardPaths} for \(\alpha =()\)).

\begin{proposition}
  Suppose that \(\alpha\) is all-ones and has \(r>0\) parts. Then
  the polynomial \(D_k^\alpha(t)\) is 1 for \(k=r\), and for \(k>r\)
  this polynomial has 
  \begin{itemize}
  \item degree \(k-r\),
  \item  constant term \(1\),
  \item  leading coefficient \((-1)^{k-r+1} (k-1)!/(r-1)!\).
  \end{itemize}
\end{proposition}

As in \eqref{eq:anka} we have that 
 when \(\alpha\) is a composition
consisting of \(r\) ones,
\begin{equation}
  \label{eq:ank}
  a_{n,k}^{\alpha} \sim 
   \frac{1}{k^k} D_k^\alpha(1/k) \frac{k^{k-r}}{(k-r)!} k^n =
   \frac{D_k^\alpha(1/k)}{(k-r)!} k^{n-r}
  \qquad \text{ as } n \to \infty 
\end{equation}

\begin{remark}
  This was stated  incorrectly in 
  \cite[Corollary 4]{Snellman:StandardPaths}; 
  the numerator was evaluated at 1 rather
  than at \(1/k\).
\end{remark}

\begin{remark}
We have not been able to determine a formulae for the value of 
\(D_k^\alpha(1/k)\).
\end{remark}

Although the \emph{poles} of the rational function \(L_k\) is of
greater interest than the zeroes (since the pole of smallest modulus,
namely \(1/k\), determines the asymptotic growth of the Taylor
coefficients), we could still ask where the zeroes are located. By 
\eqref{eq:Lksa}, the zeroes of \(L_k\) are \(0\) together with the
zeroes of \(D_k\). We make the following conjecture:

\begin{conjecture}
  There is some \(\R \ni c \approx 8\) and a curve \(C \subset \C\)
  such that, when \(k\) is large, the zeroes of \(D_k^{(1)}(x/k)\) are
  either close to 
  the set \(\setsuchas{k/m}{k \in \Nat^+} \cap [c,k]\) or lie
  interspersed close to the curve \(C\).

  Thus the zeroes of \(D_k^{(1)}(x)\) are either close to
  \(\setsuchas{1/m}{1 \in \Nat^+} \cap [c/k,1]\) or lie intersperesed
  close to the curve \(k^{-1}C\).
\end{conjecture}
The zeroes of \(D_k^{(1)}(x/k)\) is shown in Figure~\ref{fig:limcyc},
and those zeroes that approach the curve \(C\) is shown in greater
detail in Figure~\ref{fig:limcyc2}.
\begin{figure}[t]
  \centering
  \includegraphics[bb=117 150 494 700, height=7truecm,
  width=9truecm]{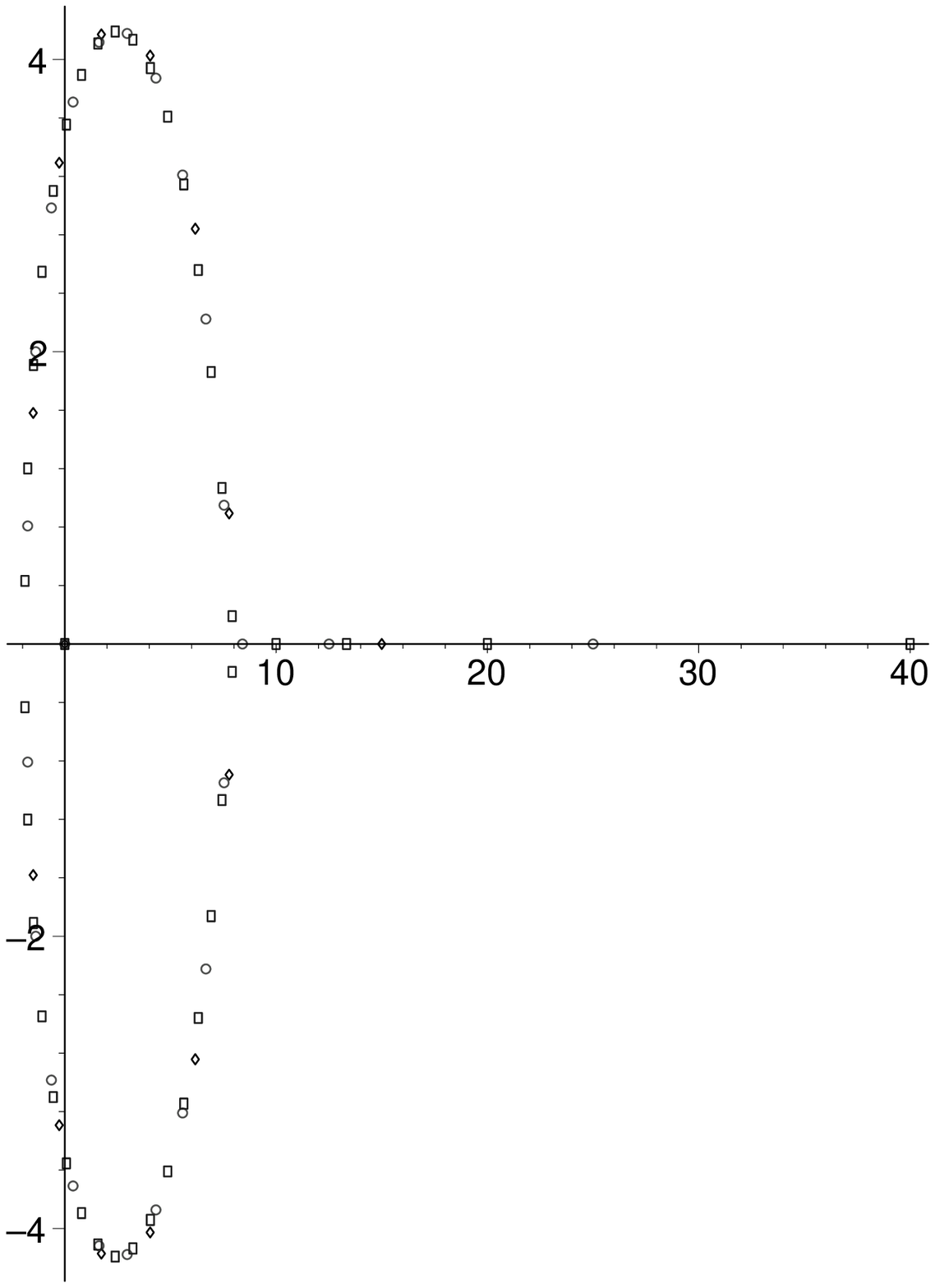} 
  \caption{The zeroes of \(D_k^{(1)}(x/k)\) for \(k=15,25,40\).}
  \label{fig:limcyc}
\end{figure}

\begin{figure}[t]
  \centering
  \includegraphics[bb=80 80 530 710, height=8truecm,
  width=8truecm,clip]{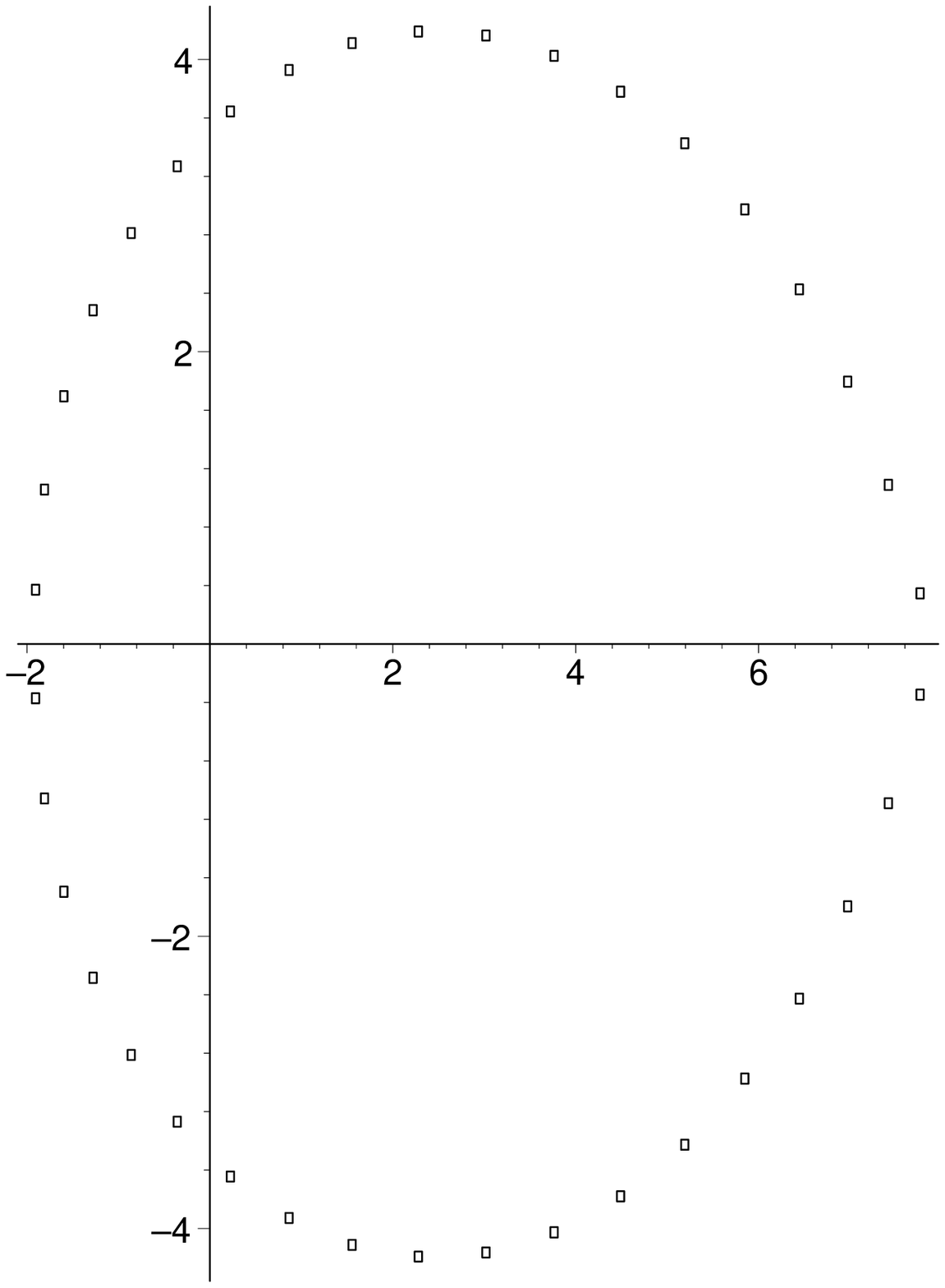} 
  \caption{The non-sporadic zeroes of \(D_{45}^{(1)}(x/k)\).}
  \label{fig:limcyc2}
\end{figure}

    \end{subsection}

    \begin{subsection}{Enumeration of saturated chains of fixed width
        in the posets \protect\(\protect\STA^d\protect\)} 
      The generating functions \(f_k^{()}[\BBD] = f_k^{(1)}[\BBD]\)
      were studied in \cite{StPa}. The authors derived an explicit
      formula for the coefficient 
      \begin{displaymath}
        [x_1^{a_1} \cdots x_k^{a_k}] f_k^{(1)}[\BBD](x_1,\dots,x_k)
      \end{displaymath}

      Recall that \(\BBD = \STA^2\) and that \(\STA^d\) is an
      increasing family of posets on \(\Compositions\), with union
      \(\STA^\infty\). Since the generating functions
      \(f_k^{\alpha}[\STA^d]\), for \(d < \infty\), satisfies the
      recurrence relation \eqref{eq:starec}, we'll be able to give
      some simple results about these functions. We note for instance
      that \(f_k^{\alpha}[\STA^d]\) are rational functions for \(d <
      \infty\). 
      Furthermore, for fixed \(k,\alpha\) it holds that 
      \begin{displaymath}
        \lim_{d \to \infty} f_k^{\alpha}[\STA^d]  = f_k^{\alpha}[\STA^\infty]
      \end{displaymath}
      in the natural formal topology on \(\Z[[x_1,\dots,x_k]]\).

 Similarly to
Theorem~\ref{thm:fkalpha} one can show 

\begin{theorem}\label{thm:stadden}
Let \(\alpha\) be a composition with \(r\) parts, and let \(k,d\) be
positive integers, such that \(r \le k\). Then the following hold:
\begin{enumerate}
\item  The denominator of \(f_k^{\alpha}[\STA^d]\)
is of the form  
\begin{equation}
  \label{eq:denstad}
  \prod_{\emptyset \neq S \subseteq \set{1,\dots,k}} 
  \left( 1 - \sum_{i \in S} x_i \right)^{e(\alpha,k,d,S)} 
\end{equation}
where \(e(\alpha,k,d,S)\) are non-negative integers, with
\[e(\alpha,k,d,\set{1,2,\dots,k}) = 1.\]
\item The denominator of
\begin{displaymath}
  L_k^{\alpha}[\STA^d](t) = f_k^{\alpha}[\STA^d](t,\dots,t)
\end{displaymath}
is of the form
\begin{equation}
  \label{eq:denta}
  \prod_{i=1}^k \left( 1 - it\right)^{c(\alpha,k,d,i)} 
\end{equation}
with \(c(\alpha,k,d,k)=1\).
\item The coefficient of \(t^n\) in \(L_k^{\alpha}[\STA^d](t)\),
  i.e. the number of saturated chains of length \(n\), starting from
  \(\alpha\), grow as (some constant times) \(k^n\) with \(n\).
\end{enumerate}
\end{theorem}

\begin{example}
Let us look at some small examples, for \(\alpha=(1)\).
We have that 
\[f_1[\STA^d] = {\frac {x_{{1}}}{1-x_{{1}}}} \]
for all \(d \ge 2\), hence that
\[f_1[\STA^\infty] = {\frac {x_{{1}}}{1-x_{{1}}}}. \]
Furthermore,
      \begin{equation}
        \label{eq:stadex}
        \begin{split}
          f_1[\STA^\infty] &= {\frac {x_{{1}}}{1-x_{{1}}}} \\
          f_2[\STA^2] &=          {\frac {x_{{1}}x_{{2}}\left( 
                1-x_{{1}}x_{{2}}\right)}{\left (1 - x_{{1}} -x_{{2}}\right) 
              \left( 1 -x_{{2}}\right) \left(  1 -x_{{1}}\right )}} \\
          f_2[\STA^3] &=
          {\frac {x_{{2}}x_{{1}}\left(
                1 -x_{{1}}{x_{{2}}}^{2}\right )}{\left(
                1 -x_{{1}}-x_{{2}}\right) \left( 1 -x_{{2}}\right) 
              \left (1 -x_{{1}}\right )}} \\
          f_2[\STA^4] &= {\frac {x_{{2}}x_{{1}}\left
                (1 -x_{{1}}{x_{{2}}}^{3}\right) }{\left(
                1 -x_{{1}} -x_{{2}}\right )\left (1 -x_{{2}}\right )\left
                ( 1 -x_{{1}} \right )}} \\
          f_2[\STA^\infty] &={\frac {x_{{2}}x_{{1}}}{\left (1-x_{{1}}+x_{{2}}
              \right)\left( 1-x_{{2}}\right) \left(1-x_{{1}}\right )}}
        \end{split}
      \end{equation}

\end{example}

From the above example, it might look like \(f_k[\STA^\infty]\) should
be rational for all 
\(k\). This is in fact not the case. Already for \(k=3\) the
denominators fail to stabilize:
the numbers \(e(1,3,d,S)\) are shown in Table~\ref{tab:STA3}.

\begin{table}[t]
  \centering
\begin{tabular}{|c|ccccccc|}
  \hline
  \(d\) & 1 & 2 & 12 & 3 & 13 & 23 & 123 \\
  2 & 2& 1& 1& 2& 1& 1& 1 \\
  3 & 3& 1& 1& 3& 1& 1& 1 \\
  4 & 4& 1& 1& 4& 1& 1& 1 \\
  5 & 5& 1& 1& 5& 1& 1& 1 \\ \hline
\end{tabular}
  \caption{the numbers \(e(1,3,d,S)\)}
  \label{tab:STA3}
\end{table}

We see that for large \(d\) the denominator of \(f_3[\STA^d]\) is of the form
\begin{multline*}
\left (1-x_{{2}}\right )
\left (1-x_{{1}}-x_{{2}}\right )
\left(1 -x_{{1}} -x_{{3}}\right )
\left (1 -x_{{2}} -x_{{3}}\right )  \times \\
\left (1 -x_{{1}} -x_{{2}} -x_{{3}}\right )
\left (1 -x_{{1}}\right )^{d}
\left (1 -x_{{3}}\right )^{d}
\end{multline*}
This means that 
\[f_3[\STA^\infty] = \lim_{d \to \infty} f_3[\STA^d]\] is not a
rational function. Similarly,  the specialization
\(L_3[\STA^\infty]\) is not rational, since \(L_3[\STA^d]\) has a
denominator of the form \((1-t)^{d+1}(1-2t)(1-3t)\).

\begin{subsubsection}{The poset \protect\(\protect\BBD\protect\)}

For  \(d=2\), i.e. for the \(\BBD\) poset, the numbers
\(e(\alpha,k,2,S)\) , for \(k=3,4\), are shown in Table~\ref{tab:BBD3}.
and Table~\ref{tab:BBD4}.
\begin{table}[t]
  \centering
\begin{tabular}{|c|ccccccc|}
  \hline
  \(\alpha\) & 1 & 2 & 12 & 3 & 13 & 23 & 123 \\
  (1) & 2& 1& 1& 2& 1& 1& 1 \\
  (2) & 2& 1& 1& 2& 1& 1& 1 \\
  (3) & 2& 1& 1& 2& 1& 1& 1 \\
  (1,1) & 1& 0& 1& 1& 1& 1& 1 \\
  (2,1) & 1& 0& 1& 0& 1& 1& 1 \\
  (2,2) & 0& 0& 1& 0& 1& 1& 1 \\
  (3,2) & 0& 0& 1& 0& 1& 1& 1 \\
  (4,4) & 0& 0& 1& 0& 1& 1& 1 \\
  (1,1,1) &0& 0& 0& 0& 0& 0& 1 \\   \hline
\end{tabular}
  \caption{the numbers \(e(\alpha,3,2,S)\)}
  \label{tab:BBD3}
\end{table}

The general pattern seems to be quite involved, even if we concentrate on
\(\alpha=(1)\), i.e. on standard paths. However, in \cite{StPa} an
explicit, though intricate
formula for the coefficient of \(x_1^{a_1}\cdots x_k^{a_k}\) in
\(f_k^{()}[\BBD]=  f_k^{(1)}[\BBD]\) is given.

\begin{table*}[t]
  \centering
  \begin{tabular}{|c|cccc cccc cccc ccc|}
    \hline
    \(\alpha\) & 1 & 2 & 3 & 4 & 5 & 6 &
7 & 8 & 9 & 10 & 11 & 12 & 13 & 14 & 15
\\ \hline
(1) & 2& 2& 2& 2& 1& 1& 1& 2& 2& 1& 1& 2& 1& 1& 1 \\
(2) & 2& 2& 2& 2& 1& 1& 1& 2& 2& 1& 1& 2& 1& 1& 1 \\
(3) & 2& 2& 2& 2& 1& 1& 1& 2& 2& 1& 1& 2& 1& 1& 1 \\
(1,1) & 1& 1& 2& 1& 1& 1& 1& 1& 2& 1& 1& 2& 1& 1& 1 \\
(1,2) & 0& 0& 2& 1& 1& 1& 1& 1& 2& 1& 1& 2& 1& 1& 1 \\
(1,3) & 0& 0& 2& 1& 1& 1& 1& 1& 2& 1& 1& 2& 1& 1& 1 \\
(1,4) & 0& 0& 2& 1& 1& 1& 1& 1& 2& 1& 1& 2& 1& 1& 1 \\
(2,2) & 0& 0& 2& 0& 1& 1& 1& 0& 2& 1& 1& 2& 1& 1& 1 \\
(2,3) & 0& 0& 2& 0& 1& 1& 1& 0& 2& 1& 1& 2& 1& 1& 1 \\
(5,9) & 0& 0& 2& 0& 1& 1& 1& 0& 2& 1& 1& 2& 1& 1& 1 \\
(1,1,1) & 0& 0& 1& 0& 0& 0& 1& 0& 1& 0& 1& 1& 1& 1& 1 \\
(1,1,2) & 0& 0& 0& 0& 0& 0& 1& 0& 1& 0& 1& 1& 1& 1& 1 \\
(1,2,2) & 0& 0& 0& 0& 0& 0& 1& 0& 0& 0& 1& 1& 1& 1& 1 \\
(5,5,5) & 0& 0& 0& 0& 0& 0& 1& 0& 0& 0& 1& 0& 1& 1& 1 \\
(1,1,1,1) & 0& 0& 0& 0& 0& 0& 0& 0& 0& 0& 0& 0& 0& 0& 1 \\
\hline
  \end{tabular}
  \caption{The numbers \(e(\alpha,4,2,S)\). Here, the columns,
    corresponding to \(\emptyset \neq S \subset{1,2,3,4}\), are coded
    in binary, e.g. \(5= 2^2+2^0\) correspond to \(\set{1,3}\).}
  \label{tab:BBD4}
\end{table*}

The specializations \(L_k^{(1)}[\BBD](t)\) looks like follows:
\begin{equation}
  \label{eq:bbdspec}
  \begin{split}
  L_1^{(1)} &= \frac{t}{1-t} \\
  L_2^{(1)} &= {\frac {\left (t+1\right ) {t}^{2}}{\left (1-t\right )
      \left (1 -2\,t \right )}} \\
  L_3^{(1)} & = {\frac {{t}^{3}\left( 3\,{t}^{2}-4\,t-1\right)}{
      \left( 1 -t\right)^{2} \left (1-2\,t \right) \left
        (1 -3\,t\right )}} \\ 
  L_4^{(1)} &= {\frac {\left
        (12\,{t}^{4}-19\,{t}^{3}-19\,{t}^{2}+13\,t+1\right ) 
      {t}^{4}}{ \left(1-4\,t\right) \left(1 -t\right)^{2} 
      \left(1 -2\,t\right)^{2} \left(1-3\,t\right)}} \\
  \end{split}
\end{equation}

Unfortunately, the recurrence relation \eqref{eq:bbdrec} for
\(f_k[\BBD]\) does not
specialize to a recurrence relation for \(L_k[\BBD]\) in the way
that the recurrence relation \eqref{eq:ncrec} for \(f_k[\NC]\) does, so
even if one should be able to guess the general form of \(L_k[\BBD]\)
it would not be trivial to prove it.

\end{subsubsection}
    \end{subsection}

    \begin{subsection}{Enumeration of shadow skew tableaux in
        \protect\(\protect\NC\protect\)  and
        \protect\(\protect\BBD\protect\)} 
      Let \(\mathcal{Q}\) denote either the poset \(\NC\) or the poset
      \(\BBD\). For these two \(\Young\)-graded posets, we have
      defined (in subsection \ref{sec:shadow}) the shadow of a tableau
      encoding a saturated chain: this is a skew tableau encoding a
      saturated chain in the Young lattice. Conversely, for a
      saturated chain in the Young lattice, we have defined its
      multiplicity as the number of saturated chains in
      \(\mathcal{Q}\) having the \(\Young\)-chain as its shadow.

      Let \(\Sort: \C[[x_1,\dots,x_k]] \to \C[[x_1,\dots,x_k]]\) be
      the continuous, \(\C\)-linear map defined on monomials by
      \[\Sort(\boldsymbol{x}^{\boldsymbol{\alpha}}) =
      \boldsymbol{x}^{\boldsymbol{\beta}},\]
      where \(\boldsymbol{\beta} \in \Nat^k\) is the \emph{dominant
        weight} associated to \(\boldsymbol{\alpha}\), i.e. the
      entries in \(\boldsymbol{\alpha}\) sorted in decreasing
      order.

      Now suppose that \(\boldsymbol{\beta} \in \Young\) has all parts
      equal, so that there is only one composition which has the same
      parts. Then it is clear that the generating functions for
      saturated chains in \Young, starting from \(\boldsymbol{\beta}\),
      and counted with multiplicity \(m(\gamma)\), is given
      by
      \begin{equation}
        \label{eq:cmul}
        \begin{split}
        \tilde{f}_k^{\beta}[\mathcal{Q}](x_1,\dots,x_k) &:=
      \sum_{\substack{\gamma 
          \text{ saturated chain in {\Young } of width } k \\ 
          \text{ starting from } \boldsymbol{\beta} }}
      m(\gamma)v(\gamma) \\
      &= 
      \Sort\bigl(f_k^{\beta}[\mathcal{Q}](x_1,\dots,x_k)\bigr)          
        \end{split}
    \end{equation}

    \begin{example}
      Let \(\mathcal{Q}=\NC\), \(\alpha=()\). Then 
      \begin{multline}
        \label{eq:tild}
        \tilde{f}_2 = \Sort(f_2) = \Sort\Bigl(
          {\frac{x_{{1}}x_{{2}}\left(
                1 -x_{{1}}x_{{2}}\right )}{\left (1 -
                x_{{1}}\right )\left (1 -x_{{2}}\right)
              \left(1 -x_{{1}} -x_{{2}}\right)}}
        \Bigr) 
          \\ =
         \Sort\Bigl(
         x_{{1}}x_{{2}} + 2\,{x_{{1}}}^{2}x_{{2}} +
         2\,x_{{1}}{x_{{2}}}^{2} + 3\,{x_{{1}}}^{3}x_{{2}} +
         4\,{x_{{1}}}^{2}{x_{{2}}}^{2} + 3\,x_{{1}}{x_{{2}}}^{3}
         + \dots \Bigr)  \\ =
         x_1x_2 + x_1^2\left(4\,x_{{2}}+4\,{x_{{2}}}^{2}\right) +
         x_1^3\left(6\,x_{{2}}+14\,{x_{{2}}}^{3}+14\,{x_{{2}}}^{2}\right)
         + \cdots
      \end{multline}
      That the coefficient of \(x_1^2x_2\) is 4 is consistent with the
      fact that there are 4 standard paths in \(\NC\) that ends either
      in \((2,1)\) or in \((1,2)\), namely the standard paths which
      has diagrams shown in Table~\ref{tab:sha4}.

        \begin{table}[t]
          \centering
        \begin{tabular}{cccc}
          \setlength{\unitlength}{0.5cm}
          \begin{picture}(3,5)
            \multiput(0,0)(0,1){2}{\PBox}
            \multiput(1,0)(0,1){1}{\PBox}
            \put(0.2,0.2){1}
            \put(0.2,1.2){3}
            \put(1.2,0.2){2}
          \end{picture}    
        &
          \setlength{\unitlength}{0.5cm}
          \begin{picture}(3,5)
            \multiput(0,0)(0,1){2}{\PBox}
            \multiput(1,0)(0,1){1}{\PBox}
            \put(0.2,0.2){1}
            \put(0.2,1.2){2}
            \put(1.2,0.2){3}
          \end{picture}    
          &
          \setlength{\unitlength}{0.5cm}
          \begin{picture}(3,5)
            \multiput(0,0)(0,1){1}{\PBox}
            \multiput(1,0)(0,1){2}{\PBox}
            \put(0.2,0.2){2}
            \put(1.2,1.2){3}
            \put(1.2,0.2){1}
          \end{picture}    
          &
          \setlength{\unitlength}{0.5cm}
          \begin{picture}(3,5)
            \multiput(0,0)(0,1){1}{\PBox}
            \multiput(1,0)(0,1){2}{\PBox}
            \put(0.2,0.2){3}
            \put(1.2,1.2){2}
            \put(1.2,0.2){1}
          \end{picture}    
        \end{tabular}
          \caption{Standard paths ending in a composition with shadow
            \((2,1)\).} 
          \label{tab:sha4}
        \end{table}
    \end{example}
      
    We conjecture that the series \(\tilde{f}_k^\alpha[\mathcal{Q}]\)
    are non-rational in all non-degenerate cases.
      
    \end{subsection}
\end{section}


\begin{section}{Labeled enumeration of saturated chains of fixed
    width}
  \begin{subsection}{Labeling the edges of the Hasse diagram}
    Let \(\mathcal{Q}\) be one of the posets \(\NC\), \(\BBD\) or
    \(\STA^d\). 
  We label the edges in the Hasse diagram of \(\mathcal{Q}\) with
  \(L,R,U_j,V_i^j\), 
  according to the type of covering relation. Saturated chains are
  labeled with the sequence of labels occurring along the edges.
  Admissible words for a composition
  \(\alpha\) now correspond bijectively to saturated chains starting from
  \(\alpha\).  In  Figure~\ref{fig:HasseNC4l} we show the the labeling of
  the edges of \(\NC\).

  \setlength{\unitlength}{1cm}
  \begin{figure*}[t]
    \begin{center}
      \begin{picture}(12,6) 
        \put(5.6,1.2){\(L\)}
        \put(6,1){\line(0,1){1}}
        
        \put(5.3,2){\(L\)}
        \put(6,2){\line(-2,1){2}}
        \put(6,2){\line(2,1){2}}

        \put(3.0,2.7){\(L\)}
        \put(4,3){\line(-3,1){3}}
        \put(4,3){\line(1,1){1}}
        \put(4,3){\line(3,1){3}}

        \put(7.5,2.1){\(U_1\)}
        \put(8,3){\line(2,1){2}}
        \put(8,3){\line(-1,1){1}}
        \put(8,3){\line(-3,1){3}}

        \put(0.2, 4.1){\(L\)}
        \put(1,4){\line(-1,1){1}}
        \put(1,4){\line(1,1){1}}
        \put(1,4){\line(3,1){3}}
        \put(1,4){\line(5,1){5}}

        \put(4.1,3.9){\(U_2\)}
        \put(5,4){\line(-1,1){1}}
        \put(5,4){\line(1,1){1}}
        \put(5,4){\line(2,1){2}}
        \put(5,4){\line(3,1){3}}

        \put(7,4){\line(-5,1){5}}
        \put(7,4){\line(-3,1){3}}
        \put(7,4){\line(0,1){1}}
        \put(7,4){\line(3,1){3}}

        \put(10.7,4.1){\(U_1\)}
        \put(10,4){\line(-2,1){2}}
        \put(10,4){\line(0,1){1}}
        \put(10,4){\line(2,1){2}}

        \put(8.5,3.1){\(U_1\)}
        \put(7.9,3.4){\(R\)}
        \put(6.4,3.2){\(L\)}




        
        


      \end{picture}

      \caption{The Hasse diagram of \(\NC\). Edges are labeled
        according to the type of the covering relation}
      \label{fig:HasseNC4l}
    \end{center}
  \end{figure*}
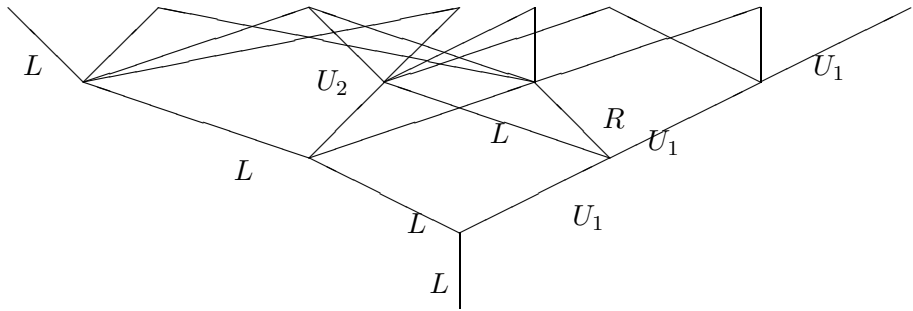

  If \(\alpha \in \Compositions\) and \(W=W_rW_{r-1}\dots W_1\) is a word 
  which is admissible for \(\alpha\), then  we give the corresponding 
  chain \(\gamma=(\alpha,W_1.\alpha, W_2W_1.\alpha,\dots, W.\alpha)\)
  non-commutative 
  weight
  \begin{equation}
    \label{eq:V}
    V(\gamma) = v(\gamma)W 
  \end{equation}
  were \(v(\gamma)\) is as in \eqref{eq:vwidth}. The non-commutative
  generalization of \eqref{eq:fk} is 
  \begin{equation}
    \label{eq:Fk}
    F_k^{\alpha}[\mathcal{Q}] = F_k^{\alpha} = \sum_{\gamma}    V(\gamma) 
  \end{equation}
  were the sum is over all saturated chains \(\gamma\) of width
  \(k\) that starts from \(\alpha\).
  Here, the \(x_i\)'s commute with each other and with the variables
  \(R,L,U_j\), but the latter variables do not commute with each
  other. Note that \(F_k^{\alpha}\) only involves finitely many variables.

    One observes that the coefficient in \(F_k^{\alpha}\) of a
    non-commutative monomial \(W\) is a single monomial in
    \(x_1,\dots,x_k\), namely the monomial encoding the endpoint of
    the path encoded by \(W\). Similarly
    the coefficient in \(F_k^{\alpha}\) of a
    commutative monomial \(\mathbf{x}^{\mathbf{a}}\) is a
    non-commutative polynomial in
    \(L,R,U_j,V_i^j\) with non-negative coefficients, encoding all
    paths (from the starting composition) that ends at
    \({\mathbf{a}}\).
     As an example, for
    \(\mathcal{Q} = \NC\) the
    coefficient of \(x_1x_2\) in \(F_2^{()}\) is \(U_2L^2 + LU_1L\).

    Clearly, specializing all non-commutative variables in
    \(F_k^\alpha\) to one
    gives \(f_k^\alpha\). On the other hand, specializing all
    commutative variables to one gives a formal power series in
    non-commuting variables, all whose occuring coefficients are
    one. If \(\mathcal{F}\) denotes the free monoid on the relevant
    non-commuting variables, then this power series is the generating
    function of the \emph{language} 
    \begin{equation}
      \label{eq:1}
      \langle \mathcal{F}; \,\, \alpha \rangle \subset \mathcal{F}.
    \end{equation}

\end{subsection}

    \begin{subsection}{Labeled enumeration in the poset
        \protect\(\protect\NC\protect\)}
      The poset \(\NC\) has covering relations given by the partial
      action of the free monoid \((\AlphabetLR \cup \AlphabetU)^*\).
      The generating function \(F_k^\alpha[\NC] = F_k^\alpha\) has
      commuting variables \(x_1,\dots,x_k\) and non-commuting
      variables in \(\set{L,R} \cup \AlphabetU\). In fact, no \(U_j\)
      with \(j > k\) 
      will occur in \(F_k^\alpha\), hence we regard \(F_k^\alpha\) as
      having non-commuting variables in
      \begin{equation}
        \label{eq:4}
        \set{L,R,U_1,U_2,\dots,U_k}
      \end{equation}
      
    \begin{theorem}\label{thm:labenumwi}
      The non-commutative generating function for labeled saturated
      chains in \(\NC\), starting from the composition
      \(\alpha=(a_1,\dots,a_s)\), 
      satisfies the recurrence
      \begin{equation}
        \label{eq:labrec}
        \begin{split}
        F_k^\alpha &= F_k^\alpha(x_1,\dots,x_k; L,R,U_1,U_2,\dots,U_k) \\
        &=
  \begin{cases}
          0 & \text{ if } k < s \\
          A +   v(\alpha) & \text{ if } k = s \\
          A + B + C
          & \text{ if } k > s \text{ and } \alpha \text{ not all-ones}
          \\
          A + B + C -  D 
          & \text{ if } k > s \text{ and } \alpha \text{ all-ones}
        \end{cases}
        \end{split}
      \end{equation}

      where 
      \begin{equation}
        \begin{split}
           A &= (x_1U_1 + \cdots + x_kU_k) F_k^\alpha \\ 
           B &= x_1L \cdot \UpSubs(F_{k-1}^\alpha)      \\
           C &= x_kR \cdot F_{k-1}^\alpha \\ 
           D & = RL^{k-1} v(\alpha)
        \end{split}
        \end{equation}
    \end{theorem}
    \begin{proof}
      This follows from Lemma~\ref{lemma:therecurs} in the same way
      that Lemma~\ref{lemma:ncrec} follows.
    \end{proof}

    For \(\alpha= (2)\), we get that 
    \begin{equation}
      \label{eq:Ftva}
      \begin{split}
        F_0^{(2)} &= 0 \\
        F_1^{(2)} &= (1-x_1U_1)^{-1}x_1^2 \\
        F_2^{(2)} &= (1-x_1U_1 -x_2U_2)^{-1}  \times \\
        & \quad \left[
          x_1L(1-x_2U_1)^{-1}x_2^2 + x_2R(1-x_1U_1)^{-1}x_1^2 \right] 
      \end{split}
    \end{equation}
    
    It is known that non-commutative rational series in finitely many
    variables  are    \emph{recognizable}, so that the coefficients
    correspond to the labels of walks from a start node to an end node
    in a certain labeled digraph. As an example, 
    \begin{multline}
        F_1^{(2)} = (1-x_1U_1)^{-1}x_1^2 = x_1^2 + x_1^3U_1 +
        x_1^4U_1^2 + \cdots
    \end{multline}
    corresponds to paths from  \(\bullet\) to  \(\circ\) in
    the following digraph: 
    \begin{displaymath}
      \color{blue}
      \xymatrix{
        \\
        \circ \ar@(ul,ur)[]^{x_1U_1} \\
        \bullet \ar[u]_{x_1^2}
        }
    \end{displaymath}

    An immediate consequence of \eqref{eq:labrec} is the following:

    \begin{theorem}\label{thm:ncdigraph}
      Let \(\alpha\) be a composition with \(r\) parts.
      Denote the language defined by \(F_k^\alpha(1,\dots,1,L,R,U_1,\dots,U_K)\)
      by \(\mathcal{L}_k^\alpha\). This is a regular language, and
    abusing notation by equating a regular language to some regular
    expression that defines it, we can write  

      \begin{equation}
        \label{eq:5}
        \mathcal{L}_r^\alpha =   (U_1 + U_2 + \dots + U_r)^*
      \end{equation}

      \begin{enumerate}[(A)]
      \item If \(\alpha\) is not all-ones, then
      a digraph for
      \(F_k^\alpha\), which enumerates 
      saturated chains of widht \(k\) in \(\NC\), starting from
      \(\alpha\),  by walks from \(\bullet\) to \(\circ\), is obtained from the
      one for \(F_{k-1}^\alpha\) by
      
      \begin{equation}
        \color{blue}
         \xymatrix{ 
          & \circ 
          \ar@(l,ul)[]^{x_1U_1} 
          \ar@(ul,u)[]^{x_2U_2} 
          \ar@(ur,r)[]^{x_kU_k} 
          & 
          \\
          *++[F]{F_{k-1}^\alpha} \ar [ur]^{x_kR}
          & 
          &  *++[F]{\UpSubs(F_{k-1}^\alpha)} \ar [ul]_{x_1L}
          \\
          & \bullet \ar [ul] \ar [ur]& \\
        }
      \end{equation}      
      Here,  
      \begin{math}
        \xymatrix{
          *++[F]{F_{k-1}^\alpha} 
          }
      \end{math}
      denotes the digraph yields \(F_{k-1}^\alpha\), and
      \begin{math}
        \xymatrix{
          *++[F]{\UpSubs(F_{k-1}^\alpha)}
          }
      \end{math}
      denotes the digraph which yields
      \(\UpSubs(F_{k-1}^\alpha)\); this digraph is  obtained
      from the former by transforming each label using \(\UpSubs\).

      It follows that 

      \begin{equation}
        \label{eq:nono}
        \mathcal{L}_k^\alpha =   
        (U_1 + U_2 + \dots + U_k)^* (L + R)  \mathcal{L}_{k-1}^{\alpha} \qquad
        \text{ for \(k \ge r\)}     
      \end{equation}

    so that
      \begin{multline}
        \label{eq:7}
        \mathcal{L}_k = \left( \sum_{i=1}^k U_i \right)^* 
        (L + R)
        \left( \sum_{i=1}^{k-1} U_i \right)^* 
        (L + R) \cdots \\ \cdots
        \left( \sum_{i=1}^{r+1} U_i \right)^* 
        (L + R)
        \left( \sum_{i=1}^{r} U_i \right)^* 
      \end{multline}

    \item If \(\alpha\) is  all-ones, then
      a digraph for
      \(F_k^\alpha\), which enumerates 
      saturated chains of widht \(k\) in \(\NC\), starting from
      \(\alpha\),  by walks from \(\bullet\) to \(\circ\), is obtained from the
      one for \(F_{k-1}^\alpha\) by
      \begin{equation}
        \color{blue} \xymatrix{
          & \circ 
          \ar@(l,ul)[]^{x_1U_1} 
          \ar@(ul,u)[]^{x_2U_2} 
          \ar@(ur,r)[]^{x_kU_k} 
          & 
          \\
          \\
          \circ 
          \ar@/^2pc/ [uur]^{x_1U_1} 
          \ar [uur]^{x_2U_2} 
          \ar@/_2pc/ [uur]_{x_kU_k}
          \\
          *++[F]{F_{k-1}^\alpha} \ar [u]^{x_kR}
          & 
          & *++[F]{\UpSubs(F_{k-1}^\alpha)} \ar [uuul]_{x_1L}
          \\
          & \bullet \ar [ul] \ar [ur]& \\
        }
      \end{equation}      
      
      hence,  for \(k \ge r\),
      \begin{equation}
        \label{eq:55}
        \begin{split}
          \mathcal{L}_k^\alpha &=   
          (U_1 + U_2 + \dots + U_k)^* \left(L + (U_1 + U_2 \dots +
            U_k) R \right) 
          \mathcal{L}_{k-1}^{\alpha} \\
          &=
          (U_1 + U_2 + \dots + U_k)^* (L + R)
          \mathcal{L}_{k-1}^{\alpha} - RL^{k-1}
        \end{split}
      \end{equation}

      \end{enumerate}
    \end{theorem}

    \begin{example}
      Since 
      \begin{displaymath}
        \xymatrix{
          *++[F]{F_{1}^{(2)}}
          }
          = 
      \color{blue} \xymatrix{
        \circ \ar@(ul,ur)[]^{x_1U_1} \\
        \bullet \ar[u]_{x_1^2}
        }
      \end{displaymath}
      and 
      \begin{displaymath}
        \xymatrix{
          *++[F]{\UpSubs(F_{1}^{(2)})}
        }
          = 
      \color{blue} \xymatrix{
        \circ \ar@(ul,ur)[]^{x_2U_1} \\
        \bullet \ar[u]_{x_2^2}
        },
      \end{displaymath}
    a digraph for  \(F_2^{(2)}\) is
    \begin{displaymath}
        \xymatrix{
          *++[F]{F_{2}^{(2)}}
          }
          = 
          \color{blue} \xymatrix{
         & \circ \ar@(l,u)[]^{x_1U_1} \ar@(r,u)[]_{x_2U_2} \\
        \cdot \ar@(dl,ul)[]^{x_1U_1}  \ar[ru]^{x_2R} 
        & 
        & \cdot  \ar@(dr,ur)[]_{x_2U_1} \ar[lu]_{x_1L}\\
        & \bullet \ar[lu]^{x_1^2} \ar[ru]_{x_2^2}
        }
    \end{displaymath}

    and 
    \begin{equation}
      \label{eq:6}
      \mathcal{L}_2^{(2)} =  (U_1 + U_2)^* (L +R)  U_1^* 
    \end{equation}

Note that specializing \(x_1=x_2=1\) in \eqref{eq:Ftva} gives
\begin{multline}
  (1-U_1 -U_2)^{-1} \left[L(1-U_1)^{-1} + R(1-U_1)^{-1} \right] = \\
  =
(1-U_1 -U_2)^{-1} (L+R)(1-U_1)^{-1}     
  \end{multline}
which correspond  precisely to \eqref{eq:6}.
    \end{example}

    \begin{example}
      A digraph for
      \(F_3^{(2)}\) is 
      \begin{displaymath}
        \color{blue} \xymatrix{
          \\
         & \cdot \ar@(l,u)[]^{x_1U_1} \ar@(r,u)[]_{x_2U_2} 
         \ar @/^3pc/ [ddrrr]^{x_3R} && \\
        \cdot \ar@(dl,ul)[]^{x_1U_1}  \ar[ru]^{x_2R} 
        & 
        & \cdot  \ar@(dr,ur)[]_{x_2U_1} \ar[lu]_{x_1L} &&\\
        & \bullet \ar[lu]^{x_1^2} \ar[ru]_{x_2^2}
        \ar[ld]^{x_2^2} \ar[rd]_{x_3^2}
         &&& \circ  \ar@(u,ur)[]^{x_1U_1} 
         \ar@(ur,r)[]^{x_2U_2} 
         \ar@(r,dr)[]^{x_3U_3} 
         \\
        \cdot \ar@(dl,ul)[]^{x_2U_1}  \ar[rd]_{x_3R} 
        & 
        & \cdot  \ar@(dr,ur)[]_{x_3U_1} \ar[ld]^{x_2L} && \\
         & \cdot \ar@(d,l)[]^{x_2U_1} \ar@(d,r)[]_{x_3U_2} 
         \ar @/_3pc/ [uurrr]_{x_1L} &&\\
        \\
        }
    \end{displaymath}
      and 
      \begin{equation}\label{eq:reglang}
        \mathcal{L}_3^{(2)} = (U_1 + U_2 + U_3)^* 
        ( L + R) 
        (U_1 + U_2)^* 
        (L +R) 
        U_1^* 
      \end{equation}
    \end{example}

\end{subsection}

\begin{subsection}{Labeled enumeration in the poset
    \protect\(\protect\BBD\protect\)} 
  In the poset \(\BBD\), we label the edges of the Hasse diagram with
  \(L\), \(U_j\) or \(V_r^1=V_r\). Defining
  \(F_k^\alpha=F_k^\alpha[\BBD]\) as before, we get:

    \begin{theorem}\label{thm:labenumwa}
      The non-commutative generating function for labeled saturated
      chains in \(\BBD\), starting from the composition
      \(\alpha=(a_1,\dots,a_s)\), 
      satisfies the recurrence
      \begin{equation}
        \label{eq:labreca}
        \begin{split}
        F_k^\alpha &= F_k^\alpha[\BBD](x_1,\dots,x_k;
        L,U_1,U_2,\dots,U_k,V_2,\dots,V_k) \\ 
        &=
  \begin{cases}
          0 & \text{ if } k < s \\
          A +   v(\alpha) & \text{ if } k = s \\
          A + B + \sum_{i=2}^k C_i
          & \text{ if } k > s 
        \end{cases}
        \end{split}
      \end{equation}

      where 
      \begin{equation}
        \begin{split}
           A &= (x_1U_1 + \cdots + x_kU_k) F_k^\alpha \\ 
           B &= x_1L \UpSubs(F_{k-1}^\alpha) \\
           C_i & = x_{i} V_i \UpSubs_i\left(F_{k-1}-
          \Delta_{i-1}^1(F_{k-1}^\alpha)\right) 
        \end{split}
        \end{equation}
    \end{theorem}
    
    \begin{example}
      We have that
      \begin{equation}
        \label{eq:qaa}
        \begin{split}
          F_1^{(1)} &= x_1(1-x_1U_1)^{-1}\\
          &= x_1 + x_1^2U_1 + x_1^3U_1^2 + x_1^4U_1^3 + \dots \\
          F_2^{(1)} &= \left[1-x_1U_1 -x_2U_2\right]^{-1}
          \Biggl(
          x_1 \left[ x_2(1-x_2U_1)^{-1} \right] \\
          & \qquad +x_2 V_2\left( x_1^2 U_1(1-x_1U_1)^{-1}\right)
          \Biggr)
        \end{split}
      \end{equation}
      
    \end{example}

    \begin{theorem}
      Let \(\alpha\) be a composition with \(r\) parts.
      Put
      \[F_k^\alpha=F_k^\alpha[\BBD](1,\dots,1,L,V_2,
      \dots,V_k,U_1,\dots,U_k).\]  
      Then a digraph for
      \(F_k^\alpha\), which enumerates 
      saturated chains of widht \(k\) in \(\NC\), starting from
      \(\alpha\),  by walks from \(\bullet\) to \(\circ\), is obtained from the
      one for \(F_{k-1}^\alpha\) by

      \begin{equation}
        \color{blue} \xymatrix{
          &&&& \circ 
          \ar@(l,ul)[]^{x_1U_1} 
          \ar@(ul,u)[]^{x_2U_2} 
          \ar@(ur,r)[]^{x_kU_k} 
          & 
          \\
          \\
          & *++[F]{B} \ar [uurrr]^{x_1L}
          & *++[F]{C_2} \ar [uurr]_{x_2V_2}  
          &
          &
           \ar@{.>} [uu]
          &
          &
          *++[F]{C_k} \ar [uull]_{x_kV_k} 
          \\
          \\
          &&&& \bullet \ar [uulll] \ar [uull] \ar@{.>} [uu]  \ar[uurr]& \\
        }
      \end{equation}      
      where \(\xymatrix{*++[F]{B}}\) is the digraph for
      \(\UpSubs(F_{k-1}^\alpha)\) and  
      \(\xymatrix{*++[F]{C_i}}\) is the digraph for 
      \begin{equation}
        \label{eq:BBB}
         \UpSubs_i
         \left(
           F_{k-1}^\alpha - \Delta_{i-1}^1(F_{k-1}^\alpha)
         \right)
      \end{equation}
      



    \end{theorem}

    This theorem is less informative than Theorem~\ref{thm:ncdigraph},
    since it is somewhat complicated to construct the digraph 
    generating \(  F_{k-1}^\alpha - \Delta_{i-1}^1(F_{k-1}^\alpha)\)
    given the digraph generating \(F_{k-1}^\alpha\). The process  may require the
    addition of extra edges, and is somewhat irregular.

    \begin{example}
      We have that
      \begin{displaymath}
        \xymatrix{
          *++[F]{F_{1}^{(1)}}
          }
          = 
      \color{blue} \xymatrix{
        \circ \ar@(ul,ur)[]^{x_1U_1} \\
        \bullet \ar[u]_{x_1}
        }
      \end{displaymath}
      and that
      \begin{displaymath}
        \xymatrix{
          *++[F]{F_{1}^{(1)} - \Delta_1^1(F_{1}^{(1)})}
        }
          = 
      \color{blue} \xymatrix{
        \circ \ar@(ul,ur)[]^{x_2U_1} \\
        \bullet \ar[u]_{x_1^2U_1}
        },
      \end{displaymath}
      so a digraph for  \(F_2^{(1)}\) is
    \begin{displaymath}
        \xymatrix{
          *++[F]{F_{2}^{(1)}}
          }
          = 
          \color{blue} \xymatrix{
         & \circ \ar@(l,u)[]^{x_1U_1} \ar@(r,u)[]_{x_2U_2} \\
        \cdot \ar@(dl,ul)[]^{x_1U_1}  \ar[ru]^{x_2V_2} 
        & 
        & \cdot   \ar@(dr,ur)[]_{x_2U_2} \ar[lu]_{x_1L}\\
        & \bullet \ar[lu]^{x_1^2U_1} \ar[ru]_{x_2U_1}
        }
    \end{displaymath}
    We furthermore see that 
    \begin{displaymath}
        \xymatrix{
          *++[F]{F_{2}^{(1)} - \Delta_1^1(F_{2}^{(1)})}
          }
          = 
          \color{blue} \xymatrix{
         & \circ \ar@(l,u)[]^{x_1U_1} \ar@(r,u)[]_{x_2U_2} \\
        \cdot \ar@(dl,ul)[]^{x_1U_1}  \ar[ru]^{x_2V_2} 
        & 
        & \cdot   \ar@(dr,ur)[]_{x_2U_2} \ar[lu]_{x_1U_1L}\\
        & \bullet \ar[lu]^{x_1^2U_1} \ar[ru]_{x_2U_1}
        }
    \end{displaymath}
    but that
    \begin{displaymath}
        \xymatrix{
          *++[F]{F_{2}^{(1)} - \Delta_2^1(F_{2}^{(1)})}
          }
          = 
          \color{blue} \xymatrix{
         & \circ \ar@(l,u)[]^{x_1U_1} \ar@(r,u)[]_{x_2U_2} \\ \\
        \cdot \ar@(dl,ul)[]^{x_1U_1}  \ar@/^2pc/[ruu]^{x_2^2U_2V_2} 
        & 
        & \cdot   \ar@(dr,ur)[]_{x_2U_2} 
        \ar@/^1pc/[luu]^{x_1x_2LU_1}
        \ar@/_1pc/[luu]_{x_1x_2U_2L}
        \\
        & \bullet \ar[lu]^{x_1^2U_1} \ar[ru]_{x_2U_1}
        }
    \end{displaymath}
    \end{example}
\end{subsection}
\end{section}

\bibliographystyle{plain}
\raggedright
\bibliography{journals,articles,snellman}

\end{document}